\documentclass[preprint,12pt]{elsarticle}
\newcommand\SHORTTITLE{Two models of delamination}
\def\version{6.10.2014 - revised 26.11.2015}\def\users{world}
\usepackage{amsmath,amssymb,amsthm,mathrsfs,epsf,a4,color,graphicx,eucal}
\usepackage{subfigure}
\usepackage{psfrag}
\usepackage{import}
\newcounter{myfigure}
\newenvironment{my-picture}[3]{\refstepcounter{myfigure}\label{#3}\setlength{\unitlength}{\textwidth}\begin{picture}(#1,#2)}{\end{picture}}

\usepackage{ifthen}\usepackage{fancyhdr}
\ifthenelse{\equal{\users}{world}}{}
{\pagestyle{fancy}\headheight=28pt
\definecolor{grey}{rgb}{0.6,0.6,0.6}
\rhead{\color{grey}\SHORTTITLE\\by C.G.Panagiotopoulos \& T.Roubicek \& V. Mantic}
\chead{} \lhead{\color{grey}Version \version, file:
\jobname.tex (in dir. ``trsevilla9'')\\compiled:
\number\day.\number\month.\number\year\ at
\the\hour:\ifnum\minute<10 0\fi\the\minute\ h }

\newcount\hour \newcount\minute
\hour=\time \divide \hour by 60 \minute=\time \loop \ifnum \minute >
59 \advance \minute by -60 \repeat
}
\usepackage[normalem]{ulem}
\definecolor{brown}{rgb}{0.5,0,0}
\ifthenelse{\equal{\users}{world}}{
    \newcommand{\REPLACE}[2]{#2}
    \newcommand{\INSERT}[1]{#1}
    \newcommand{\COMMENT}[1]{}
    \newcommand{\DELETE}[1]{}
    \newcommand{\CHECK}[1]{#1}
    
    \newcommand{\MARGINOTE}[1]{\marginpar{\color{red}\tiny\texttt{}}}
    \newcommand{\REM}[1]{\marginpar{\bfseries\tiny }}
}
{\newcommand{\REPLACE}[2]{{\color{brown}\sout{#1}\uline{#2}\color{black}}}
 \newcommand{\INSERT}[1]{{\color{blue}\uuline{#1}\color{black}}}
 \newcommand{\DELETE}[1]{{\color{brown}\sout{#1}\color{black}}}
 \newcommand{\COMMENT}[1]{{\color{red}\uuline{#1}\color{black}}}
 \newcommand{\CHECK}[1]{\color{brown}\uwave{#1}\color{black}}
 
 \newcommand{\MARGINOTE}[1]{\marginpar{\color{red}\tiny\texttt{#1}}}
 \newcommand{\REM}[1]{\marginpar{\bfseries\tiny }}

\usepackage[notcite,notref]{showkeys}
}
\definecolor{ddmagenta}{rgb}{0.7,0,0.9}
\definecolor{ddcyan}{rgb}{0,0.2,1.0}
\definecolor{dred}{rgb}{.8,0,0}

\newtheorem{theorem}{Theorem}[section]

\newtheorem{remark}[theorem]{Remark}

\numberwithin{equation}{section}

\newcommand{\TTT}{\color{black}}
\newcommand{\EEE}{\color{black}}

\newcommand{\REDD}[1]{{\color{black}#1}}

\newcommand{\DDD}[3]{\begin{array}[t]{c}#1\vspace*{-1em}\\_{#2}\vspace*{-.5em}\\_{#3}\end{array}}
\newcommand{\ddd}[3]{\DDD{\begin{array}[t]{c}\underbrace{#1}\vspace*{.6em}\end{array}}{\text{\footnotesize #2}}{\text{\footnotesize #3}}}

\newcommand\R{\mathbb R}
\newcommand\N{\mathbb N}

\newcommand{\eq}[1]{(\ref{#1})}
\newcommand\DT[1]{\mathchoice
                 {{\buildrel{\hspace*{.1em}\text{\LARGE.}}\over{#1}}}
                 {{\buildrel{\hspace*{.1em}\text{\Large.}}\over{#1}}}
                 {{\buildrel{\hspace*{.1em}\text{\large.}}\over{#1}}}
                 {{\buildrel{\hspace*{.1em}\text{\large.}}\over{#1}}}}

\newcommand\JUMP[2]{\mathchoice
                   {\big[\hspace*{-.3em}\big[#1\big]\hspace*{-.3em}\big]_{#2}}
                   {[\hspace*{-.15em}[#1]\hspace*{-.15em}]_{#2}}
                   {[\![#1]\!]_{#2}}
                   {[\![#1]\!]_{#2}}}
\newcommand\JMP[1]{\mathchoice
                   {\big[\hspace*{-.3em}\big[#1\big]\hspace*{-.3em}\big]}
                   {[\hspace*{-.15em}[#1]\hspace*{-.15em}]}
                   {[\![#1]\!]}
                   {[\![#1]\!]}}
\newcommand\bbA{\mathbb A}
\newcommand\bbC{\mathbb C}
\newcommand\bbD{\mathbb D}

\newcommand\bbT{\mathbb T}
\renewcommand\d{\mathrm d}

\newcommand\overlineGC{\overline{\Gamma}_{\mbox{\tiny\rm C}}}
\newcommand\overlineSC{\overline{\Sigma}_{\mbox{\tiny\rm C}}}

\newcommand{\fRM}{f}

\newcommand{\calD}{\mathcal{R}}

\newcommand{\calE}{\mathcal{E}}

\newcommand{\aein}{\text{a.e. in}}

\newcommand{\GD}{\mathchoice
                  {\Gamma_{\hspace*{-.15em}\mbox{\tiny\rm D}}}
                  {\Gamma_{\hspace*{-.15em}\mbox{\tiny\rm D}}}
                  {\Gamma_{\hspace*{-.1em}\mbox{\tiny\rm D}}}
                  {\Gamma_{\hspace*{-.05em}\mbox{\tiny\rm D}}}}
\newcommand{\GN}{\Gamma_{\hspace*{-.15em}\mbox{\tiny\rm N}}}
\newcommand\GC{\mathchoice{\Gamma_{\hspace*{-.15em}\mbox{\tiny\rm C}}}
                          {\Gamma_{\hspace*{-.15em}\mbox{\tiny\rm C}}}
                          {\Gamma_{\hspace*{-.05em}\mbox{\tiny\rm C}}}
                          {\Gamma_{\hspace*{-.05em}\mbox{\tiny\rm C}}}}
\newcommand{\Sdir}{\Sigma_{\mbox{\tiny\rm D}}}
\newcommand{\Snew}{\Sigma_{\mbox{\tiny\rm N}}}
\newcommand\SC{\Sigma_{\mbox{\tiny\rm C}}}
\newcommand{\psiG}{\psi_{\scriptscriptstyle\rm G}}
\newcommand{\dela}{\bbA}

\newcommand{\pwc}[2]{\overline{#1}_{\kern-1pt#2}}
\newcommand{\upwc}[2]{\underline{#1}_{\kern-1pt#2}}

\newcommand{\pwl}[2]{#1_{\kern-1pt#2}}

\newcommand{\Omegaone}{\Omega_+}
\newcommand{\Omegatwo}{\Omega_-}

\newcommand{\norm}{\vec{n}_{\mbox{\tiny\rm C}}^{}}
\newcommand{\normm}{\vec{\nu}_{\mbox{\tiny\rm C}}^{}}

\newcommand{\calR}{\mathcal{R}}

\newcommand\wD{w_{\text{\tiny\rm D}}}
\newcommand\fN{f_{\text{\tiny\rm N}}}

\begin{document}


\noindent{\LARGE\bf
\REPLACE{Two models 
of mixed mode delamination of quasistatic adhesive contact 
}
{Two adhesive-contact models for quasistatic mixed-mode delamination problems}}

\bigskip
\bigskip\bigskip

\noindent{\large\sc Christos G.\ Panagiotopoulos}\\
{\it
Institute of Applied and Computational Mathematics, Foundation for Research
and Technology - Hellas, Nikolaou Plastira 100, Vassilika Vouton, GR-700 13
Heraklion, Crete, Greece}

\medskip

\noindent{\large\sc Vladislav Manti\v c}\\
{\it Group of Elasticity and Strength of Materials,
Dept.~of Continuum Mechanics, School of Engineering,
University of Seville, Camino de los Descubrimientos s/n, ES-41092 Sevilla, Spain}

\medskip

\noindent{\large\sc Tom\'a\v s Roub\'\i\v cek}\\
{\it Mathematical Institute, Charles University, Sokolovsk\'a 83,
CZ--186~75~Praha~8,
} and
{\it Institute of Thermomechanics, Czech Academy of Sciences, Dolej\v skova~5,
CZ-182~08 Praha 8, Czech Republic}

\medskip



\bigskip

\begin{center}\begin{minipage}[t]{16.5cm}

\baselineskip=12pt
{\small \noindent {\it Abstract.}
Two models for quasistatic adhesive unilateral contact
delaminating in mixed fracture mode, i.e.\ distinguishing the
less-dissipative Mode I (opening) from the more-dissipative Mode II
(shearing), and allowing rigorous mathematical and numerical analysis,
are studied. One model, referred to as Associative Plasticity-based
Rate-Independent Model (APRIM), works for purely elastic bodies
and \DELETE{is thus rate-independent. It}\COMMENT{``rate-independent''
was repeated -- abstract should be short, not an introduction} 
involves, in addition to an 
interface damage
variable, an auxiliary variable (representing interfacial plastic slip)
to provide a fracture-mode sensitivity.
\DELETE{. Nevertheless, this model}
\REDD{It relies on a particular} concept of force-driven local solutions (given by either
vanishing-viscosity concept or maximum-dissipation principle).
\COMMENT{the reason for this concept explained in the paper - again this
is an abstract}
The other model, referred to as Linear Elastic - (perfectly) Brittle
Interface Model (LEBIM), 
\REDD{works}
 visco-elastic bodies 
\REDD{and} rely on a
conventional concept of weak solution and needs no auxiliary
interfacial variable. This model is directly related to a
usual phenomenological model of interface \REDD{fracture }
by Hutchinson and
Suo used in engineering. 
\REDD{This paper devises} a way how the phenomenology of the LEBIM can be fit 
to imitate the APRIM under relatively very slow loading, where both models are
essentially rate-independent.
The so-called effective dissipated energy is partitioned in both
formulations to the surface energy and the energy dissipated during the
interface debonding process, where the former is independent and the latter
dependent on the fracture mode mixity.
A numerical comparison of these models, implemented in a Boundary Element
Method (BEM) code, is
carried out on a suitable two-dimensional example. Furthermore, the
computational efficiency and behaviour of the LEBIM is illustrated on
another geometrically more complicated numerical example.


\medskip

\noindent {\it Key Words.} Inelastic surface damage,
associative rate-independent model, interfacial gradient plasticity
with hardening and damage, maximally-dissipative solution,
non-associative model, weak solution, semi-implicit discretization, visco-elastic material, quadratic mathematical programming.


\medskip

\noindent {\it AMS Subject Classification:}
35K85,  
65K15, 
65M38, 
74M15, 
74R20. 
}
\end{minipage}
\end{center}


\bigskip

\baselineskip=16pt

\section{Introduction -- two models}\label{sec-intro}

Adhesive contact represents
a part of {\it nonlinear contact mechanics}
with  numerous practical applications. \CHECK{Modelling damage and fracture of the adhesive layer}
standardly involves an ``internal'' variable distributed at
the contact surface.   \CHECK{Here, this damage-type variable is  denoted by $z$.} 
This concept, falling in the
framework of generalized standard materials \cite{HalNgu75MSG}, was devised
by M.\,Fr\'emond \cite{Fre82,Fre87}.

An important feature  appearing in  engineering models
(and so far mostly omitted in the mathematical literature), is  the
dependence of this process  on the fracture modes under which it proceeds.
Indeed, Mode I (=\,opening) usually dissipates much less energy than
Mode II (=\,shearing). The difference may be up to hundreds of percents, cf.\
\cite{BanAsk00NFCI,HutSuo92MMCL,LieCha92ASIF,Mant08DRLM,SwLiLo99ITAM}. Moreover,
the delamination process seldom occurs in such pure modes and, in reality, the
{\it  mixed fracture-mode} appears more frequently. The substantial difference
in the dissipation in different modes can be explained either by some
roughness of the glued interface (to be overcome in Mode II but not in
Mode I, cf.\ \cite{ERDC90FEBI}) or by some plastification caused by shear
in Mode II (but \REDD{much less } 
by mere tension in Mode I) before the delamination itself
happens, cf.~\cite{LieCha92ASIF,TveHut93IPMM}. Let us emphasize that
a strong dependence of dissipation (or fracture toughness) on the fracture
mode mixity (easily ranging hundreds of percents)
has been observed in extensive experiments
\cite{BanAsk00NFCI,ERDC90FEBI,HutSuo92MMCL,LieCha92ASIF,SwLiLo99ITAM} and thus
this phenomenon cannot be neglected in models.

The most common modelling assumption is that the time scale of external load
is rather slow so that, in particular the viscous bulk or surface
effects are negligible and inelastic processes  on the contact
surface are much faster so that they are considered rate-independent.
Often, also inertial effects are neglected. We also confine ourselves
to isothermal models at small strains and frictionless adhesive
contact undergoing a unidirectional delamination (i.e.\ no healing
is considered). Most
continuum-mechanical models
in literature can be (re)formulated as generalized gradient flow or,
equivalently as the doubly-nonlinear Biot-type equation:
\begin{align}\label{Biot}
\partial_{\DT q}\mathcal{R}(q;\DT q)+\partial_q\mathcal{E}(t,q)\ni0 ,
\end{align}
where $q=q(t)$ is a state and $\DT q$ denotes its time derivative,
$\mathcal{R}$ is a dissipation potential
and $\mathcal{E}$ is a stored-energy potential, and ``$\partial$''
is a partial generalized derivative (here the convex subdifferential).


Two models developed in the literature that allow for rigorous mathematical 
and numerical analysis \REDD{(as far as existence of solutions and 
numerically stable, implementable, and convergent approximation schemes)}  
are studied here. Both have a bit
different character:
\begin{itemize}
\item[LEBIM:\hspace*{-1.5em}]\hspace*{1.5em}
{\it Linear Elastic - (perfectly) Brittle Interface Model}
(like \cite{Bennati2009,cornetti12A1,MTBGP14ALEB,TMGCP10ACTA,TMGP11BEMA,VC06EFM}):
\item[$+$]\vspace*{-.7em}
It allows for incorporation of arbitrary {\it phenomenological} dependence
of dissipation (or activation) energy on the fracture mode mixity. The state $q$
in \eqref{Biot} is a couple $(u,z)$ and the ``elasticity domain'' of the
adhesive varies depending on the state, which is sometimes
called ``{\it non-associative}'', cf.\ e.g.\ \cite{HanRed99PMTN}.
Mathematic\REDD{s} 
relies on \REDD{the} conventional concept of weak solution,
cf.\ \cite{KruPaRo??MMS}.
\vspace*{-.7em}\item[$-$]
On the other hand, mathematic justification of the merely elastic
variant is not available in literature, and a modification
by considering a {\it visco-elastic} rheology of the delaminating
bodies seems inevitable, cf.\ \cite{KruPaRo??MMS}. Thus the viscosity is an
additional parameter of the model, which may enrich the modelling aspects
but, if the desired model is rather purely elastic, may a bit twist the
model.

\bigskip

\item[APRIM:\hspace*{-1.8em}]\hspace*{1.8em}
{\it Associative Plasticity-based Rate-Independent Model}
(devised in \cite{tr+mk+jz,RoMaPa13QMMD}):
\item[$+$]\vspace*{-.7em} Arising from a rational motivation of
interfacial {\it plasticity with hardening}, it uses
$\calR$ independent of $q$ and degree-1 homogeneous in $\DT q$
and works for purely
elastic bodies and elasto-plastic adhesive. It involves an
auxiliary variable $\pi$ (called \CHECK{interface plastic slip}) to execute
the \REDD{fracture-}mode-sensitivity. The state $q$
in \eqref{Biot} is a triple $(u,z,\pi)$ and the ``elasticity domain'' of the
internal variables $(z,\pi)$ of the adhesive is independent of the state,
and such attribute is sometimes denoted by the adjective ``associative''.
\vspace*{-.7em}\item[$-$]It needs \REDD{a reasonable choice of a solution
concept, as discussed for this particular model in \cite{VoMaRo14EMDL}.
To avoid unwanted effects leading to rather unphysically too early delamination
if the energy in the elastic bulk dominates, one must realize }
a special concept of force-driven local solutions (here either 
a vanishing-viscosity concept like \cite{Roub13ACVE,RoPaMA13QACV} or 
\REDD{solutions complying at least approximately with the} maximum-dissipation 
principle \cite{Roub??MDLS,RoPaMa14LSAQ,VoMaRo14EMDL}).
\end{itemize}
It should be emphasized that some other models (mostly representing
a variant of the LEBIM or the so-called Cohesive Zone Models (CZM)) are routinely used in engineering
even though any rigorous convergence/existence analysis
is not at disposal. Such computations and the models themselves are thus
unjustified so far from the mathematical point of view, although in particular cases the launched
computational simulations may be physically relevant; see e.g.\
\cite{CaDaMo03NSMM,CoCa11ES,MTBGP14ALEB,TMGCP10ACTA,TMGP11BEMA,VC06EFM}
and references therein.

Let us still mention that the LEBIM was rigorously analyzed already in
\cite{rr+tr2} even in a full thermodynamical context but exploiting the
concept of non-simple materials (see e.g.~\cite{toupin}) which whould be
much more demanding to be implemented computationally.

In this paper, after more detailed formulation of the models
LEBIM and APRIM in Sections~\ref{sec-model-A} and \ref{sec-model-B},
respectively, we devise a way how the phenomenology of the LEBIM can be fit
to imitate the APRIM under very slow loading, where both models are
essentially rate-independent  and have a good chance to produce
similar responses in spite of rather different essence of both models;
this is performed in Section~\ref{sec-A-B}. Eventually, in
Section~\ref{sec-simul},   these models, implemented in a  BEM code,  are numerically tested and compared
  on a two-dimensional  single-domain  example. Furthermore,   the solution of a two-domain problem,  reproducing the Mixed  Mode Flexure (MMF) test, by LEBIM is studied.

Let us first introduce the notation we will use
\TTT throughout \EEE the paper
which will be common to both mentioned models.
We suppose that the 
(visco)elastic\REDD{-bulk}/inelastic-adhesive
structure occupies
a bounded Lipschitz domain $\Omega\subset\R^d$ composed by
(for notational simplicity only) two bodies, denoted by $\Omegaone$ and
$\Omegatwo$, glued together on a common contact boundary, denoted by $\GC$,
which represents a prescribed
interface where delamination may occur.
This means we consider
\[
\Omega = \Omegaone \cup \GC \cup \Omegatwo\,,
\]
with $\Omegaone$ and $\Omegatwo$ \REDD{being}
disjoint Lipschitz subdomains.
We  denote by $\vec{n}$ the
outward unit normal to $\partial \Omega$, and by $\norm$
 the unit normal to $\GC$, which we consider
oriented from $\Omegaone$ to $\Omegatwo$. Moreover, given  $v$
defined on $\Omega {\setminus} \GC$,
$v^+$ (respectively, $v^-$)
signifies the restriction of $v$ to $\Omegaone$ (to $\Omegatwo$,
resp.).
 We further suppose  that the boundary of $\Omega$ splits as
\[
\partial \Omega = \GD\cup \GN\,,
\]
with $\GD$ and $\GN $ open subsets in the relative topology of
$\partial\Omega$,  disjoint one from each other, \CHECK{and for $d=3$ each of them
with a smooth (one-dimensional) boundary}. Considering $T>0$ a fixed
time horizon, we  set
\begin{displaymath}
Q:=(0,T){\times}\Omega, \quad \Sigma:= (0,T){\times}\partial
\Omega, \quad \SC\!:= (0,T){\times}\GC , \quad \Sdir\!:= (0,T){\times}
\GD, \quad \Snew\!:= (0,T){\times}\GN.
\end{displaymath}
For readers' convenience, let us summarize the basic notation used
in what follows:

\vspace{.7em}

\hspace*{-1.4em}\fbox{
\begin{minipage}[t]{
18em}
\small

$d=2,3$ dimension of the problem,

$u:\Omega{\setminus}\GC\to\R^d$ displacement,


$z:\GC\to[0,1]$ delamination variable,

$e=e(u)=\frac12\nabla u^\top\!+\frac12\nabla u$ small-strain tensor,

$\JMP{u}= \CHECK{u^-|_{\GC} - u^+|_{\GC}}\ $ jump of $u$ across $\GC$, \COMMENT{This is due to the Signorini condition for $\JMP{u} \norm\geq0$ and the above definition of $\norm$. We have discussed it several times, but the error persists;
and in particular in this manuscript it requires to change sign in several formulas, see below.}

$\sigma$  stress tensor,


$\bbC\in\R^{d^4}
$
elastic moduli in the Hook law,

$\bbA\in\R^{d\times d}$
elastic coefficients of the adhesive,

$\kappa_{\rm n}^{}$ distributed normal stiffness,

$\kappa_{\rm t}^{}$ distributed tangential stiffness,
\end{minipage}\ \
\begin{minipage}[t]{20em
}\small

$\psiG$ fracture-mode-mixity angle for the LEBIM,






$\bbD\in\R^{d^4}$ viscosity constants for the LEBIM,




$\alpha=\alpha(\JMP{u})$ effective energy dissipated on $\GC$,





\REDD{$w_{\rm D}:\Sdir\to\R^d$}
 prescribed
boundary displacement,

\REDD{$\fRM_{\rm N}:\Snew\to\R^d$ } applied traction,






$\sigma_{\mathrm{t,yield}}$ yield shear stress for the APRIM,

$\pi{:}\GC{\to}\R^{d-1}$ the interfacial plastic slip  for
APRIM,

$\kappa_{\rm H}$
modulus of kinematic hardening for APRIM,

$a_{_{\rm I}}$ energy released per unit  area in pure Mode I,\\
$a_{_{\rm II}}$ energy released per unit  area in pure Mode II

\end{minipage}\medskip
}

\vspace{-.2em}

\begin{center}
{\small\sl Table\,1.\ }
{\small
Summary of the basic notation used 
\TTT throughout \EEE the paper.
}
\end{center}

\TTT Throughout \EEE
the whole paper, we will assume  \REDD{that}
$\bbC,\,\bbD:\R_{\mathrm{sym}}^{d\times d}\to\R_{\mathrm{sym}}^{d\times d}$
and $\bbA:\R^d\to\R^d$ \REDD{are} linear positive definite, and $\alpha(\cdot)>0$.

\section{The non-associative model (LEBIM)}\label{sec-model-A}



\noindent The {\it state} is formed by the couple $(u,z)$.
We use {\it Kelvin-Voigt's rheology} \REDD{for subdomains (adherents)}  and, rather for mathematical
reasons to facilitate analysis in multidimensional cases,
a (possibly only slightly) nonlinear static response.
Hence we assume the \emph{stress}
\REDD{$\sigma:Q\rightarrow\R^{d\times d}$} in the form:
\begin{align}\label{KV-stress}
\sigma=\sigma(u,\DT{u}) :=\!\!\!\ddd{\bbD
e(\DT{u})}{viscous}{stress}\!\!\!
+\!\!\!\ddd{\bbC e(u)}{elastic}{stress}\!\!\!.
\end{align}
Furthermore, we shall denote by $T_{\rm n}$ and $T_{\rm t}$ the normal and the
tangential components of the traction $\sigma\big|_{\Gamma}\vec{n}$ on some
two-dimensional surface $\Gamma$
(used for either $\Gamma=\GC$ or $\Gamma=\GN$), i.e.
defined on $\GC\cup\GN$ respectively by the formulas
\begin{align}\label{def-of-T}
T_{\rm n}=\vec{n} \cdot\sigma\big|_{\Gamma}\vec{n}
\qquad\text{and}\qquad
T_{\rm t}=\sigma\big|_{\Gamma}\vec{n}-T_{\rm n}\vec{n}
\end{align}
where of course we take as $\vec{n}$ the unit normal  $\norm$ to $\GC$.
\TTT Note that $T_{\rm n}$ is a scalar while $T_{\rm t}$ is a vector.\EEE

The classical formulation of the adhesive contact problem
consists, beside the force equilibrium \REDD{(neglecting body forces)} 
supplemented  with standard boundary conditions, from two complementarity 
problems on $\SC$. Altogether, we have the boundary-value problem
\begin{subequations}\label{adhes-classic}
\begin{align}
\label{eq6:adhes-class-form1} &
\mathrm{div}\big(\sigma
\big)
=0,
&&\text{in }Q{\setminus}\SC,
\\
\label{eq6:adhes-class-form2} &u=\wD &&\text{on }\Sdir,\hspace{1.2em}
\\\label{eq6:adhes-class-form3-bis}
&
\sigma
\vec{n}
=\REDD{\fN} &&\text{on }\Snew,\hspace{1.2em}
\\\label{adhes-form-d1}
& \JMP{\sigma
}\norm
=0
 &&\text{on
}\SC,\hspace{1.2em}
\\
\label{adhes-form-d2-} &
T_{\rm t}
\REPLACE{+}{-}  z\big(\bbA \JMP{u}{-}
\TTT(\norm\bbA \JMP{u}){\cdot}\norm\big)\EEE
=0
&&\text{on }\SC,\hspace{-.15em}
\\
\label{adhes-form-d2} &\JMP{u}{\cdot}\norm\ge0
\ \ \ \ \text{ and }\ \ \ \ T_{\rm n}
\REPLACE{+}{-}  z(\bbA\JMP{u}){\cdot}\norm \REPLACE{\ge}{\le} 0
&&\text{on }\SC,\hspace{1.2em}
\\
\label{adhes-form-d4} &
(\JMP{u}{\cdot}\norm)\big(T_{\rm n} \REPLACE{+}{-}  z(\bbA\JMP{u})\norm
\big)=0
&&\text{on }\SC,\hspace{1.2em}
\\\label{adhes-form-d6}
&\DT{z}\le0
\ \ \ \ \text{ and }\ \ \ \ \xi\le
\alpha_{_{\rm PLAST}}(\JMP{u}) \ \ \ \ \text{ and }\ \ \ \
\DT{z} \left(\xi - \alpha_{_{\rm PLAST}}(\JMP{u})\right) =0
&&\text{on }\SC,\hspace{1.2em}
\\\label{adhes-form-d8}
& \xi\in
\mbox{$\frac12$}\bbA\JMP{u}{\cdot}\JMP{u}-\alpha_{_{\rm ADHES}}+  N_{[0,1]}(z) &&\text{on
}\SC.\hspace{1.2em}
\end{align}
\end{subequations}
\COMMENT{with reference to (2.3f): ONCE THE ADHESIVE LAYER IS BROKEN, THE NORMAL TRACTION MUST BE NON-POSITIVE}
%

\TTT In \eqref{adhes-form-d1}, $\JMP{\sigma}$ denotes naturally the jump
of the stress tensor $\sigma$ accross $\GC$ so that \eqref{adhes-form-d1}
expresses the equilibrium of \REDD{tractions}   on $\GC$. \EEE
The complementarity problem \eqref{adhes-form-d2}--\eqref{adhes-form-d4}
describes the frictionless Signorini unilateral contact.
%
The complementarity problem~\eqref{adhes-form-d6}--\eqref{adhes-form-d8}
corresponds to the \emph{flow rule} governing the evolution of $z$:
\begin{equation}
\label{e:reactivation-2}
\frac12\bbA\JMP{u}{\cdot}\JMP{u} + N_{(-\infty,0]}(\DT{z}) + N_{[0,1]}(z)
\ni\alpha_{_{\rm ADHES}}+\alpha_{_{\rm PLAST}}(\JMP{u})=:\alpha(\JMP{u})\qquad\text{in $\SC$,}
\end{equation}
with
$N_{(-\infty,0]},\,N_{[0,1]}: \R \rightrightarrows\R$
denoting the normal cones (in the sense of convex analysis) to the
intervals $(-\infty,0]$ and $[0,1]$, respectively.
{F}or more details about derivation of the model we refer to
\cite{rr+tr2,rr+tr}. We will consider the initial-value problem for
\eqref{adhes-classic} by prescribing the initial condition
\begin{align}\label{IC}
u(0)=u_0 \quad \aein \ \Omega\qquad\text{ and }\qquad z(0)=z_0 \quad
\aein \ \GC.
\end{align}

In the weak formulation, the boundary-value problem \eqref{adhes-classic}
takes the abstract form \eqref{Biot}, i.e.\ now
the doubly-nonlinear Biot-type system of two inclusions:
\begin{align}\label{Biot-A}
\partial_{\DT u}\mathcal{R}(\DT u)+\partial_u\mathcal{E}(t,u,z)\ni0
\ \ \ \ \ \text{ and }\ \ \ \ \
\partial_{\DT z}\mathcal{R}(u;\DT z)+\partial_z\mathcal{E}(u,z)\ni0,
\end{align}
provided $\calE$ and $\calR$ are taken as
\begin{align}
&\calE(t,u,z):=\begin{cases}
\displaystyle{
\TTT\frac12\EEE
\!\int_{\Omega{\setminus}\GC}\!\!\!\!
\bbC e(u){:}e(u)
\,\d x
-\int_{\GN}\REDD{\fN}(t){\cdot}u\,\d S\!}&
\\[-.2em]\displaystyle{\hspace{2.5em}
+\frac12\!\int_{\GC}\!\!\!
z\dela\JMP{u}{\cdot}\JMP{u}-\alpha_{_{\rm ADHES}}z
\,\d S} &
\text{if }u=\wD(t)\ \text{ on }\GD,\\[-.5em]
&\text{$\JMP{u}{\cdot}\norm\ge0$
and $0\le z\le1$ on $\GC$,}
\\[-.0em] + \infty & \text{otherwise,}\end{cases}
\label{8-1-k}
\\
\textrm{and} & \nonumber
\\&
\label{DISS} \calD\big(u;\DT u,\DT z):=
\begin{cases}
\displaystyle{
\TTT\frac12\EEE
\int_{\TTT\Omega{\setminus}\GC\EEE}\!\!\!\bbD e(\DT{u}){:}e(\DT{u})\,\d x
+\int_{\GC}\!\alpha_{_{\rm PLAST}}(\JMP{u})|\DT z|\,\d S}  & \text{if
$\DT z\le0$ a.e. in $\GC$,}
\\[-.2em] + \infty & \text{otherwise.}
\end{cases}
\end{align}
\COMMENT{Is $\overlineGC$ used in the second integral of (2.8) intentionally? Has it been defined?}
Note that the {\it effective dissipated energy}
\begin{align}\label{alpha-split}
\alpha(\JMP{u})=\alpha_{_{\rm ADHES}}+\alpha_{_{\rm PLAST}}(\JMP{u})
\end{align}
is, in general, composed from
a part $\alpha_{_{\rm ADHES}}$ interpreted as an {\it adhesion energy} stored by
creating a new surface by delamination and another part, denoted
by $\alpha_{_{\rm PLAST}}$,
interpreted as an {\it energy dissipated by the plastification} process
during delamination. Altogether,
$\alpha_{_{\rm ADHES}}+\alpha_{_{\rm PLAST}}$ is to be understood as a
{\it fracture energy} (sometimes also referred to as {\it fracture toughness}).  Only the plastification energy is considered 
\REDD{as}
dependent on the delamination mode, and is usually increasing from pure Mode I
to pure Mode II, and need not be zero even in pure Mode I,
cf.\ \cite[Fig.\,4.2]{Brac12IA}, \cite[Fig.\,1]{ERDC90FEBI}  or
\cite[Fig.\,2]{VBMG03NMIF}.
As the delamination is considered here as a unidirectional process, the part
$\alpha_{_{\rm ADHES}}$ cannot be refreshed back and is effectively dissipated,
too. \CHECK{Both parts would be distinguished if $\calR$ would be finite }
(i.e.\ healing or re-bonding of the adhesive would be allowed) or if the
full thermodynamical context would be considered (then only
$\alpha_{_{\rm PLAST}}$ but not $\alpha_{_{\rm ADHES}}$ would contribute
to the heat production).


A standard engineering approach, as e.g.~in
\cite{HutSuo92MMCL,MTBGP14ALEB,TMGCP10ACTA,TMGP11BEMA}, is that
$\dela$ is diagonal \CHECK{in a local coordinate system associated to $\GC$, writing e.g. $\bbA={\rm diag}(\kappa_{\rm n}^{},\kappa_{\rm t}^{},\kappa_{\rm t}^{})$ for $\norm=(1,0,0)$ at some $x\in\GC$},
 and the
activation energy denoted here by $\alpha$,
\CHECK{ whereas in engineering typically referred to as fracture energy and denoted  by $G_c$, }
depends on the so-called \emph{energetic fracture-mode-mixity angle} $\psiG$ defined as
\begin{align}\label{a-mixity-dependence-}
\psiG=\psiG\big(\JMP{u}\big):=
{\rm arc\,tan}\bigg(
\sqrt{\frac{\kappa_{\rm t}}{\kappa_{\rm n}}}
\frac{|\JUMP{u}{\rm t}|}{|\JUMP{u}{\rm n}|}\bigg),
\end{align}
\DELETE{usually lower than one, - I DO NOT UNDERSTAND THIS}
where $\JUMP{u}{\rm t}$ and $\JUMP{u}{\rm n}$ stand for
the tangential and the normal relative displacements; i.e.\
the jump of displacement across the boundary $\GC$ decomposes
as $\JMP{u}= \JUMP{u}{\rm n}\norm+\JUMP{u}{\rm t}$, with
$\JUMP{u}{\rm n}=\JMP{u}\cdot\norm$.
\CHECK{Occasionally, other fracture-mode-mixity angles may be  defined  as associated to   displacements or
stresses,  $\psi_u={\rm arc\,tan}(|\JUMP{u}{\rm t}|/|\JUMP{u}{\rm n}|)$ and $\psi_{\sigma}={\rm arc\,tan}(\kappa_{\rm t}|\JUMP{u}{\rm t}|/\kappa_{\rm n}|\JUMP{u}{\rm n}|)$,
respectively. Here, we will not  use any of  the   phenomenological laws for $G_c(\psiG)$ well-known in engineering \REDD{(see \cite{MTBGP14ALEB})} but
rather fit it with the plasticity-inspired model described in the
following section. }

\section{\REDD{The} associative plasticity-based
model (APRIM)}\label{sec-model-B}

\def\Z{z}

In accordance with other experimental observations  and computational models
\cite{LieCha92ASIF,TveHut93IPMM,EXPPLSL2},
 the associated plastic zones in the
adjacent bulk, near the crack tip, are larger in Mode II than in Mode I and
these plastic phenomena are localized in a relatively narrow plastic zone in
the bulk in the interface vicinity. In order to provide a better representation
of these experimental results, a plastic tangential response has been assumed
at the interface, which allows us to distinguish between fracture Mode I and II
in the sense that some additional dissipated energy is associated to interface
fracture in Mode II. In our works \cite{tr+mk+jz,RoMaPa13QMMD,PaMaRo13CM},
imitating the conventional models of linearized single-threshold plasticity
with kinematic hardening (e.g.\ \cite{HanRed99PMTN}), we have invented a
plastic slip variable $\pi$ as another internal variable on the delaminating
surface (in addition to $z$) whose evolution activates in Mode II but not
in Mode I in order to \REDD{achieve the desired mixity-mode-dependent 
dissipation}.
In contrast to the LEBIM from Section~\ref{sec-model-A}, the
mathematical analysis now needs a gradient theory used for some of
the internal variables; here, following \cite{RoPaMa14LSAQ}, we apply it
on $\pi$ rather than on $z$ as in \cite{tr+mk+jz,RoMaPa13QMMD}.

Instead of \eqref{KV-stress}, we  consider purely elastic material now:
\begin{align}\label{Hook-stress}
\sigma=\sigma(u) :=\bbC e(u).
\end{align}
In the classical formulation, referring to \eqref{def-of-T},
this model uses again (\ref{adhes-classic}a-d)
completed now by  
\begin{subequations}\label{adhes-classic-APRIM}
\begin{align}
\label{adhes-form-d2-*} &
T_{\rm t} \REPLACE{+}{-} z\big(w{-}\big(w{\cdot}\norm\big)\norm\big)=0
\ \ \text{ with }\ w:=\TTT\bbA\EEE(\JMP{u}-\bbT\pi)
&&\text{on }\SC,\hspace{-.15em}\\\label{adhes-form-d2*} &
\JMP{u}{\cdot}\norm\ge0
\ \ \text{ and }\ \ T_{\rm n} \REPLACE{+}{-} z\,w{\cdot}\norm \REPLACE{\ge}{\le} 0
\ \ \text{ and }\ \ \big(\JMP{u}{\cdot}\norm\big)\big(T_{\rm n} \REPLACE{+}{-} z
\,w
{\cdot}\norm\big)=0
&&\text{on }\SC,
\\
&\DT{z}\le0\ \text{ and }\ \xi\le a_1\
\text{ and }\
\DT{z}\left(\xi{-}a_1\right) =0\ \text{ with }\
\xi\!\in\!\mbox{$\frac12$}w{\cdot}\TTT\bbA\EEE^{-1}w-a_0+N_{[0,1]}(z)\!\!
&&\text{on }\SC,
\\&\nonumber
\DT\pi=\begin{cases}
0&\text{if }|\zeta|<\sigma_{\mathrm{t,yield}}.
\\[-.2em]\lambda\zeta,\ \ \lambda\ge0&
\text{if }|\zeta|=\sigma_{\mathrm{t,yield}},
\end{cases}\ \ \ \REDD{\text{ and }\ \ |\zeta|\le\sigma_{\mathrm{t,yield}}}\ \ \text{ with }
\\
&\label{class-form4+}
\qquad\qquad\qquad\qquad\quad\ \zeta=
\REDD{T_{\rm n}}
-{\rm div}_{\mbox{\tiny S}}(\kappa_{_{\rm G}}\!\nabla_{\mbox{\tiny S}}\pi)
+({\rm div}_{\mbox{\tiny S}}\norm)(\kappa_{_{\rm G}}\!\nabla_{\mbox{\tiny S}}\pi\norm)
-\kappa_{_{\rm H}}\pi&&\text{on }\SC,
\end{align}
\COMMENT{SHOULD WE DEFINE $\TTT\bbA\EEE$?}
where $\bbT=\bbT(x):\R^{d-1}\to\R^d$ denotes
the mapping from the space where $\pi$ has values to the
tangent space of the $d$-dimensional space where
$\JMP{u}$ has values. \REDD{This mapping allows to sum up 
$(d{-}1)$-dimensional vector $\pi$ with the $d$-dimensional 
vector $\JMP{u}$, as needed in \eqref{adhes-form-d2-*} and later in 
\eqref{E-delam-small-k-II}, too.}
In \eqref{class-form4+},
$\kappa_{_{\rm G}}>0$ is a (presumably small) parameter determining length-scale
of spatial variation of $\pi$.
Also note that
$\xi$ and $\zeta$ denote the available driving \TTT forces \EEE for the
activated evolution of $\Z$ and $\pi$, respectively.
In \eqref{class-form4+},
${\rm div}_{\mbox{\tiny S}}:={\rm trace}(\nabla_{\mbox{\tiny S}})$
denotes the $(d{-}1)$-dimensional ``surface divergence'' and
$\nabla_{\mbox{\tiny S}}$ a ``surface
gradient'', i.e.\ the tangential derivative defined as
$\nabla_{\mbox{\tiny S}}v=\nabla v-(\nabla v{\cdot}\norm)\norm$ for
$v$ defined on $\GC$.

In addition to the boundary conditions (\ref{adhes-classic}b,c),
due to the surface gradient of $\pi$, we now need also the condition for
$\pi$. Most naturally, we assume:
\begin{align}
\label{class-form-9}
\nabla_{\mbox{\tiny S}}\pi{\cdot}\TTT\normm\EEE=0\ \ \ \,\text{ on }\partial\GC,
\end{align}
where
$\TTT\normm\EEE$ denotes the  unit vector lying in $\GC$ and being outward
normal to $\partial\GC$.
Of course, the boundary/transmi\TTT{ss}\EEE{ion}-value problem
(\ref{adhes-classic}a--d) with
\eqref{Hook-stress} and (\ref{adhes-classic-APRIM}a--e) is to be completed
by the initial conditions. Assuming undamaged and unplastified adhesive at
the beginning, we consider
\begin{align}\label{IC-APRIM}
\pi(0,\cdot)=\pi_0=0\qquad\text{ and }\qquad z(0,\cdot)=z_0=1.
\end{align}\end{subequations}
\TTT Note that we are now choosing a specific $z_0$ in \eqref{IC} and, \EEE
in contrast to \eqref{IC}, the initial condition for $u$ is now irrelevant
or, more precisely, it follows from \eqref{IC-APRIM} because $u(0)$ is
assumed to solve the boundary/transmition-value problem (\ref{adhes-classic}a--d)
with \eqref{Hook-stress} and (\ref{adhes-classic-APRIM}a,b).
The last term in
\eqref{class-form4+} involves $({\rm div}_{\mbox{\tiny S}}\REDD{\norm})$ 
which is (up to a factor $-\frac12$) the mean
curvature of the surface $\GC$, and it arises
by applying a Green's formula on a curved surface
\begin{align}\label{Green-surface}
\int_{\GC}\!\!\TTT\widetilde{v}\EEE{:}((\nabla_{\mbox{\tiny S}}v){\otimes}\norm)\,\d S
=\int_{\GC}\!\!({\rm div}_{\mbox{\tiny S}}\norm)(\TTT\widetilde{v}\EEE{:}(\norm{\otimes}\norm))v
-{\rm div}_{\mbox{\tiny S}}(\TTT\widetilde{v}\EEE{\cdot}\norm)v\,\d S
+\int_{\partial\GC}\!\!(\TTT\widetilde{v}\EEE{\cdot}\TTT\normm\EEE)v\,\d l
\end{align}
used with
$\TTT\widetilde{v}\EEE=
\kappa_{_{\rm G}}\!\nabla_{\mbox{\tiny S}}\pi$
and with the boundary condition $\TTT\widetilde{v}\EEE{\cdot}\TTT\normm\EEE=0$ on $\partial\GC$,
cf.\ \eqref{class-form-9}, for the directional-derivative term
$\int_{\GC}\kappa_{_{\rm G}}\!\nabla_{\mbox{\tiny S}}\pi{\cdot}\nabla_{\mbox{\tiny S}}\tilde\pi\d S$ arising from the term
$\int_{\GC}\frac12\kappa_{_{\rm G}}|\nabla_{\mbox{\tiny S}}\pi|^2\d S$
in \eqref{E-delam-small-k-II}.
The Green-type identity \eqref{Green-surface}
was used in a similar context in mechanics of complex (also called non-simple)
continua,
cf.~\cite{toupin,ppg-mv}.

This quite complicated boundary-value problem can still be covered by the
simple and elegant  Biot's form \eqref{Biot} even with $\calR(q,\DT q)$
independent of $q$, which is why we used the adjective ``associative'',
provided the state $q=(u,z,\pi)$ and provided the stored-energy functional is
taken as
\begin{subequations}\label{delam-small-II}
\begin{align}
\label{E-delam-small-k-II}
& {\calE}(t,u,z,\pi):=\!\left\{\begin{array}{ll}
\displaystyle{
\int_{\Omega{\setminus}\GC}\!
\frac12\bbC e(u){:}e(u)\,\d x
+\int_{\GC}\!\!\!
 \Big(
\TTT\frac12\EEE\!\Z\TTT\bbA\EEE\big(\JMP{u}\!{-}\bbT\pi\big){\cdot}\big(\JMP{u}\!{-}\bbT\pi\big)-a_0\Z
\!}
\hspace{-3.5em}&\\[.5em]
\displaystyle{\hspace{1.2em}+\,\frac{\kappa_{_{\rm H}}}{2}|\pi|^2
+\frac{\kappa_{_{\rm G}}}{2}\big|\nabla_{\mbox{\tiny S}}\pi\big|^2\Big)}
\displaystyle{\,\d S-\int_{\GN}\!\!\REDD{\fN}(t){\cdot}u\,\d S\!}
&\!\!\!\!\!\!\!\!\!\!\!\!\text{if }u=\wD(t)\ \text{ on }\GD,\\[-.5em]
&\!\!\!\!\!\!\!\!\!\!\!\! 0\le z\le1\text{ on }\GC, \text{ and }
\\[.3em]
&\!\!\!\!\!\!\!\!\!\!\!\!\JUMP{u}{\rm n}\ge0 \text{ on }\GC,
\\[.3em]
+\infty&\!\!\!\!\!\!\!\!\!\!\!\!\text{elsewhere,}
\end{array}\right.\hspace*{-1em}
\intertext{while the dissipation-energy
potential is taken as}
\label{R-delam-small-k-II}
&\calR(\DT z,\DT\pi):=
\begin{cases}\displaystyle{\int_{\GC}\!\!a_1\big|\DT z\big|
+\sigma_{\mathrm{t,yield}}\big|\DT\pi\big|\d S}
&\text{ if }\DT z\le0\text{ a.e.~on }\GC,\\
+\infty&\text{ otherwise}.
\end{cases}
\end{align}
\end{subequations}
Actually, the doubly-nonlinear Biot-type equation \eqref{Biot} can
now be written as three inclusions:
\begin{align}\label{Biot-B}
\partial_u\mathcal{E}(t,u,z,\pi)\!\ni\!0,
\ \ \ \
\partial_{\DT z}\mathcal{R}(\DT z)+\partial_z\mathcal{E}(u,z,\pi)\!\ni\!0
\ \, \text{ and }\ \ \ \,
\partial_{\DT\pi}\mathcal{R}(\DT\pi)+\partial_\pi\mathcal{E}(u,z,\pi)\!\ni\!0.
\end{align}

Similarly as in Section~\ref{sec-model-A},
the {\it effective dissipation energy} in Mode I, denoted here
by $a_{_{\rm I}}$, is composed from two parts, namely the {\it surface energy}
$a_0$  associated to the   creation of a new surface and the  energy dissipated
by debonding process $a_1$, i.e.
\begin{align}\label{def-of-aI}
a_{_{\rm I}}:=a_0+a_1.
\end{align}

As in \eqref{a-mixity-dependence-}, we will consider that
$\TTT\bbA\EEE={\rm diag}(\kappa_{\rm n}^{},\kappa_{\rm t}^{},\kappa_{\rm t}^{})$
\CHECK{in a local coordinate system with $\norm=(1,0,0)$. }
Starting from the initial conditions \eqref{IC-APRIM},
the response in pure Mode I is essentially determined by $\kappa_{\rm n}^{}$ and
$a_{_{\rm I}}$ because pure opening neither triggers the evolution of $\pi$ nor
causes $\JUMP{u}{\rm t}\ne0$, cf.~Fig.\,\ref{fig5:delam-mode}(Left).
To analyse the response in pure Mode II, let us realize that
the tangential stress $\sigma_{\rm t}$ is a derivative of $\calE$ with
respect to $\JUMP{u}{\rm t}$, and thus
$\sigma_t=\sigma_t(u,\pi)=\kappa_{\rm t}^{}(
\JUMP{u}{\rm t}{-}\pi)$ if $z=1$.
\COMMENT{WE SHOULD SAY SOMEWHERE THAT THE FOLLOWING IS VALID FOR $d=2$ ONLY}
In analogy with the conventional
plasticity, the slope of the evolution of $\pi$ \CHECK{vs. $\JUMP{u}{\rm t}$ }
under hardening is $\kappa_{\rm t}^{}/(\kappa_{\rm t}^{}{+}\kappa_{_{\rm H}})$.

\begin{my-picture}{.9}{.33}{fig5:delam-mode}
\psfrag{MODE-I}{\footnotesize \frame{\,MODE I$^{^{^{}}}_{}$\,}}
\psfrag{MODE-II}{\footnotesize \frame{\,MODE II$^{^{^{}}}_{}$\,}}
\psfrag{[u]n}{\footnotesize $\JUMP{u}{\rm n}$}
\psfrag{[u]t}{\footnotesize $\JUMP{u}{\rm t}$}
\psfrag{stress}{\footnotesize $\sigma_{\rm n}$}
\psfrag{stress2}{\footnotesize $\sigma_{\rm t}$}
\psfrag{sd}{\footnotesize \begin{minipage}[t]{5em}
$\!\sqrt{2\kappa_{\rm n}^{}a_{_{\rm I}}}$
\\[-.2em]$=\hspace*{-.3em}{:}\,\sigma_{\mathrm{n,crit}}$\end{minipage}}
\psfrag{sd1}{\footnotesize \begin{minipage}[t]{5em}
$\!\!\sqrt{2\kappa_{\rm n}^{}a_{_{\rm I}}}$
\\[-.2em]$=\hspace*{-.3em}{:}\,\sigma_{\mathrm{n,crit}}$\end{minipage}}
\psfrag{sd2}{\footnotesize \begin{minipage}[t]{5em}
$\!\!\sqrt{2\kappa_{\rm t}^{}a_{_{\rm I}}}$
\\[-.2em]$=\hspace*{-.3em}{:}\,\sigma_{\mathrm{t,crit}}$\end{minipage}}
\psfrag{sd3}{\footnotesize $\hspace*{-1.7em}\sigma_{\mathrm{t,yield}}$}
\psfrag{s}{\footnotesize $\sqrt{2a_{_{\rm I}}/\kappa_{\rm n}^{}}$}
\psfrag{s2}{\small $\frac{\sigma_{\mathrm{t,yield}}}{\kappa_{\rm t}^{}}$}
\psfrag{s3}{\footnotesize $u_{_{\rm II}}$ from \eqref{ch5:delam-small-s-pi-II}}
\psfrag{s4}{\footnotesize  $\sqrt{2a_{_{\rm II}}/\kappa_{\rm t}^{}\!}$}
\psfrag{area}{\footnotesize \footnotesize \begin{minipage}[t]{10em}
+\ =$\,a_{_{\rm I}}\,$=\,effective\\[-.3em]
\hspace*{4em}dissipated\\[-.4em]\hspace*{4em}energy\end{minipage}}
\psfrag{slope}{\footnotesize slope=$\kappa_{\rm n}^{}$}
\psfrag{slope1}{\footnotesize slope=$\kappa_{\rm t}^{}$}
\psfrag{slope2}{\footnotesize{slope=}{\small$\frac{\kappa_{\rm t}\kappa_{_{\rm H}}}{\kappa_{\rm t}{+}\kappa_{_{\rm H}}}$}}
\psfrag{area2}{\footnotesize area=$a_{_{\rm II}}$}
\psfrag{pi0}{\footnotesize $\kappa_{_{\rm H}}\pi_{_{\rm II}}$}
\psfrag{a0}{\footnotesize $a_0$ (stored energy)}
\psfrag{a1}{\footnotesize  $a_1$ (dissipated energy)}
\psfrag{A0}{\footnotesize $a_0$}
\psfrag{A1}{\footnotesize  $a_1$}
\psfrag{A3}{\footnotesize  \begin{minipage}[t]{15em}
$\frac12\kappa_{_{\rm H}}\pi_{_{\rm II}}^{\TTT2\EEE}=\,$energy stored\\[-.2em]
\hspace*{4.5em}via hardening\end{minipage}}
\psfrag{A4}{\footnotesize \footnotesize \begin{minipage}[t]{10em}
$\sigma_{\mathrm{t,yield}}\pi_{_{\rm II}}$=\,energy\\[-.2em]
\hspace*{2em}dissipated by\\[-.4em]\hspace*{2em}plastification\end{minipage}}
\psfrag{piII}{\footnotesize $\pi_{_{\rm II}}$}
\psfrag{evolution}{\footnotesize back-stress $\kappa_{_{\rm H}}\pi$}
\psfrag{of-pi}{\footnotesize of $\pi$}
\hspace*{-1em}\vspace*{-1.em}{\includegraphics[width=.88\textwidth]{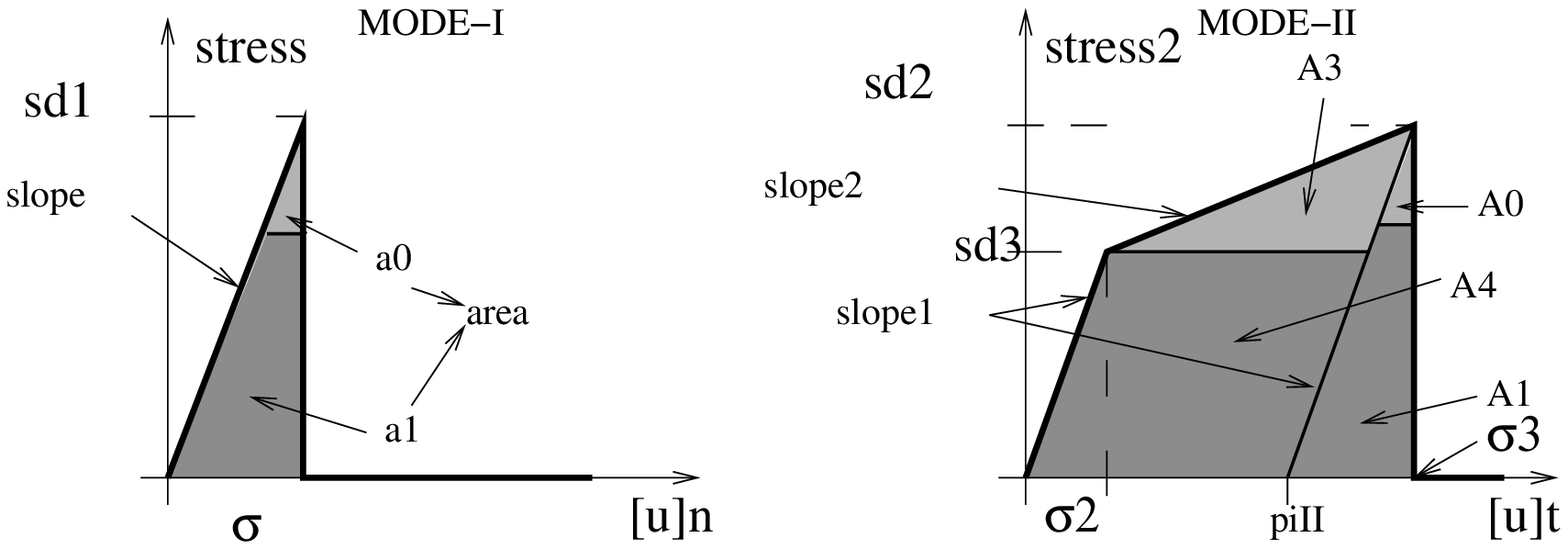}}
\end{my-picture}\nopagebreak\vspace*{-1.em}
\begin{center}
{\small Fig.\,\ref{fig5:delam-mode}.\ }
\begin{minipage}[t]{.85\textwidth}\baselineskip=8pt
{\small
Schematic illustration of the stress-relative displacement law in
the model \eqref{delam-small-II} for $d=2$ and
$\TTT\bbA\EEE={\rm diag}(\kappa_{\rm n}^{},\kappa_{\rm t}^{})$ in the
opening (Mode I) and the shearing (Mode II) experiments,
considering $z_0=1$ and $\pi_0=0$. The partition of the effective dissipated
energy is depicted only formally.
The contribution of the delamination-gradient term is neglected, i.e.\
$\kappa_{_{\rm G}}{=}0$.
}
\end{minipage}
\end{center}

{F}rom \REPLACE{\eqref{ch5:delam-small-II-activation}}{Fig.\,\ref{fig5:delam-mode}(Right)},
one can see that delamination in pure Mode II is triggered if
$\frac12\kappa_{\rm t}^{}|\JUMP{u}{\rm t}\!{-}\pi|^2
=\frac12\sigma_{\rm t}^2/\kappa_{\rm t}^{}$
attains the threshold $a_{_{\rm I}}$, i.e.~if the tangential stress
$\sigma_{\rm t}$ achieves the critical stress $\sigma_{\mathrm{t,crit}}
=\sqrt{2\kappa_{\rm t}a_{_{\rm I}}}$\DELETE{, as depicted in
Fig.\,\ref{fig5:delam-mode}(Right)}.
Delamination in Mode II is thus triggered by the tangential
displacement $u_{_{\rm II}}$ and the tangential slip $\pi_{_{\rm II}}$ given by
\begin{subequations}\label{ch5:delam-small-a-s-II}
\begin{align}
\label{ch5:delam-small-s-pi-II}
u_{_{\rm II}}&=\frac{\sqrt{2\kappa_{\rm t}a_{_{\rm I}}}
(\kappa_{\rm t}^{}{+}\kappa_{_{\rm H}})-\sigma_{\mathrm{t,yield}}\kappa_{\rm t}^{}}
{\kappa_{\rm t}^{}\kappa_{_{\rm H}}}
\ \ \ \text{ and }\ \ \
\pi_{_{_{\rm II}}}=\frac{\sqrt{2\kappa_{\rm t}a_{_{\rm I}}}
-\sigma_{\mathrm{t,yield}}}{\kappa_{_{\rm H}}}
\intertext{
and, after some algebra, one can see that the overall
 {\it effective dissipation energy} in Mode  {II}, denoted here
by $a_{_{\rm {II}}}$, is composed from four parts,
cf.~Fig.\,\ref{fig5:delam-mode}(Right), namely}
\label{ch5:delam-small-a-II}
a_{_{\rm II}}&=a_{_{\rm I}}+\sigma_{\mathrm{t,yield}}\pi_{_{\rm II}}
+\frac12\kappa_{_{\rm H}}\pi^2_{_{\rm II}}
=\ddd{a_1\,+\,\sigma_{\mathrm{t,yield}}\pi_{_{\rm II_{_{_{_{_{}}}}}}}}{arising from the}{dissipated energy $\calR$}
\!\!\!+\!\!\!
\ddd{a_0\,+\,\frac12\kappa_{_{\rm H}}\pi_{_{\rm II}}^2}{arising from the}{stored enery $\calE$}
\end{align}\end{subequations}
provided $2\kappa_{\rm t}^{}a_{_{\rm I}}\ge \sigma_{\mathrm{t,yield}}^2$, and taking
into account that the evolution of $\pi$ will stop evolving after delamination
is completed. To have actually this behaviour, the parameters should
satisfy
\begin{align}\label{2-sided-condition}
\frac12\sqrt{2\kappa_{\rm t}^{}a_{_{\rm I}}}<\sigma_{\mathrm{t,yield}}
\le\sqrt{2\kappa_{\rm t}^{}a_{_{\rm I}}}.
\end{align}
After some algebra, from \eqref{ch5:delam-small-a-s-II} it
results that $\kappa_{_{\rm H}}\pi_{_{\rm II}}=\kappa_{\rm t}^{}\kappa_{_{\rm H}}(u_{_{\rm II}}-{\sigma_{\mathrm{t,yield}}}/{\kappa_{\rm t}^{}})/(\kappa_{\rm t}^{}{+}\kappa_{_{\rm H}})$,
cf.~Fig.\,\ref{fig5:delam-mode}(Right).
A measure of maximum
\emph{fracture-mode sensitivity}\index{fracture-mode sensitivity}
$a_{_{\rm II}}/a_{_{\rm I}}$ is then indeed bigger than 1, namely
\begin{align}\label{2-conditions}
\frac{a_{_{\rm II}}}{a_{_{\rm I}}}
=1+\Big(\sigma_{\mathrm{t,yield}}\pi_{_{_{\rm II}}}+\frac{\kappa_{_{\rm H}}\pi_{_{_{\rm II}}}^2}{2}\Big)/a_{_{\rm I}}=
1+\frac{\kappa_{\rm t}}{\kappa_{_{\rm H}}} - \frac{\sigma^2_{\rm t,\rm yield}}{2\kappa_{_{\rm H}}a_{_{\rm I}}}.
\end{align}
\TTT Like discussed already for \eqref{alpha-split}, also here the particular
contributions in the spliting \eqref{def-of-aI} could be distinquished
only in a full thermodynamical model where $a_1$ (together possibly with
$\sigma_{\mathrm{t,yield}}\pi_{_{\rm II_{_{_{_{_{}}}}}}}$) would contribute to the
heat production while $a_0$ (together possibly with
$\frac12\kappa_{_{\rm H}}\pi_{_{\rm II}}^2$) would not, cf.\ also
the splitting \eqref{ch5:delam-small-a-II}.
\EEE

The vital part of the rate-independent models themselves is a suitable
choice of a concept of solutions. Let us emphasize that the class of the so-called
local (in fact weak \cite{Roub??MDLS}) solutions to nonconvex rate-independent
systems like this one given
by \eqref{delam-small-II} is typically very wide and not every type solution
will serve well in our context. For example, the so-called energetic solutions
\cite{MieThe04RIHM} which conserves energy will have a tendency for too early
rupture preferably in less dissipative mode, i.e.\ here Mode I.
Here, as we want to relate APRIM with LEBIM,
we should consider a {\it force-driven}-like {\it local solution}, which is
physically relevant and comparable with the conventional weak solution concept
for the viscous LEBIM. See \cite{Miel11altDEMF}  for a comparison of various
concepts of local solution to rate-independent systems, and  \cite{VoMaRo14EMDL} for a numerical comparison of the energetic and a particular local solution.

A physically justified concept of local solutions is a {\it vanishing-viscosity
solution}, obtained by considering small viscosity in the bulk like in the
LEBIM or/and in the adhesive and pass it to zero. The limit, if exists,
is the vanishing-viscosity solution. When combined with a suitable spatial  discretisation,
a control  of the (approximate) energy conservation needs very fine
time discretisation for small viscosities, cf.\ \cite{RoPaMA13QACV},
and therefore this method is computationally quite heavy.

Computationally very efficient but physically in general rather ad-hoc
method (cf.\ the discussion in \cite{Roub??MDLS}), which however has a
guaranteed convergence (in terms of subsequences) towards local solutions
devised in \cite{RoPaMa14LSAQ}, is based on a semi-implicit discretisation of
the fractional-step type, cf.\ \eqref{Biot-B-disc} below. \REDD{Cf.\ also 
\cite[Sect.\,4.3.4.3]{MieRou15RIST}.}

\section{Reading information from  APRIM  to \REPLACE{use for}{be used in} LEBIM
}\label{sec-A-B}
We further proceed by describing the APRIM
in terms of the \CHECK{\emph{effective dissipated energy}, also referred to as 
\REDD{a} \emph{fracture energy} $\alpha(\cdot)=G_{c}(\cdot)$},
 for which  it can
be shown \cite{MTBGP14ALEB,Krenk1992,LENCI2001,ShahinTaheri2008,Cornetti2009}
that equals the energy per unit area stored in the adhesive at the crack tip.
 In principle, one can imagine various strategies
fitting LEBIM towards APRIM,   justified from the perspective of
plastification like suggested in \cite{LieCha92ASIF,TveHut93IPMM}
  as an explanation of the fracture-mode-mixity sensitivity,
or even conversely fitting APRIM towards a given phenomenology of
LEBIM.

Here, rather as an example, we will make the fitting of LEBIM to APRIM in such a way
that Mode I will be the same for both models while, in Mode II, the
overall dissipated energy will be the
same for both models, cf.\ Fig.\,\ref{fig5:delam-mode-II}.
Two scenarios for fitting their constitutive law of interface   are considered.
In the first one, the initial  stiffness of the adhesive layer is fitted, Fig.~\ref{fig5:delam-mode-II}(Left). In the second one,
the relative tangential displacement  causing the rupture in Mode II will be equal, then the tangential
adhesive  stiffness  in the matrix $\bbA$ for LEBIM, namely
$\kappa_{\rm t}^{}=2a_{_{\rm II}}/u_{_{\rm II}}^2$ with $a_{_{\rm II}}$ and
$u_{_{\rm II}}$ from \eqref{ch5:delam-small-a-s-II}, will   be lower
than the adhesive \REDD{stiffness} $\kappa_{\rm t}$ in APRIM, Fig.~\ref{fig5:delam-mode-II}(Right).

 With reference to the first scenario, a certain misfit in the results  can be expected because the  LEBIM fitted according to Fig.~\ref{fig5:delam-mode-II}(Left) is elastically stiffer
 in Mode II than APRIM for sufficiently large stresses,  whereas the  LEBIM fitted according to Fig.~\ref{fig5:delam-mode-II}(Right) is at least initially more compliant.
 Thus, we can expect  that
 some parts in the resulting global behaviour of both models  will be better fitted in one scenario, but expectedly this will
create a certain misfit in some other parts.

One should also realize
that APRIM has some unexploited freedom: e.g.\ one can consider
a non-linear hardening, i.e.\ a non-quadratic term instead of
$\frac12\kappa_{_{\rm H}}|\pi|^2$ in \eqref{E-delam-small-k-II}. Having
used a gradient term in \eqref{E-delam-small-k-II}, this generalization
would still allow for a mathematical support. Similarly, the
elastic $\bbA$-term in \eqref{E-delam-small-k-II} can be made
non-quadratic. Obtaining such a freedom, one can then try conversely
to fit APRIM towards LEBIM with a given phenomenological nonlinearity
$\alpha(\cdot)$.

In engineering applications, the phenomenological law of $G_{\rm c}$ proposed in
\cite{HutSuo92MMCL}, cf.\ \cite{BanAsk00NFCI,MTBGP14ALEB}, is usually applied,
\begin{align}\label{a-mixity-dependence}
G_{\rm c}(\psiG)=a_{_{\rm I}}\big(1+{\rm tan}^2((1{-}\lambda)\psiG)\big),
\end{align}
where  $a_{_{\rm I}}=G_{\rm c}(0^\circ)$ gives the fracture energy in Mode I, and $\lambda$ is
the so-called mode sensitivity  parameter, $0\leq \lambda\leq 1$. A
moderately strong fracture-mode sensitivity occurs when the ratio $a_{_{\rm II}}/a_{_{\rm I}}$
   is about 5-10 (see Fig. \ref{paramratio}(a)),  with  $a_{_{\rm II}}=G_{\rm c}(90^\circ)$
the fracture energy in Mode II, which happens for $\lambda$ about 0.2-0.3.
A numerical implementation of \eqref{a-mixity-dependence} in the above non-associative
model was presented in \cite{KruPaRo??MMS}.

 In a theoretical study of the behaviour of APRIM, including an interface plastic slip
variable, the following functional dependence of $G_{\rm c}(\psiG)$  was deduced in
\cite{PaMaRo13CM,ROMAPlzni11}, while an in-depth analysis and detailed description
was presented in \cite{VoMaRo14EMDL},
\begin{align}\nonumber
\alpha(\JMP{u})&:=G_{\rm c}(\psiG(\JMP{u}))
\\&\:=
\begin{cases}\displaystyle{a_{_{\rm I}}},  &\text{for }\ 0\leq\psiG\leq\arcsin\displaystyle{\frac{\sigma_{\rm t,\rm yield}}{\sqrt{2\kappa_{t} a_{_{\rm I}}}}},\\[-.2em]
  \displaystyle{\frac{2 a_{_{\rm I}}(\kappa_{\rm t}{+}\kappa_{_{\rm H}})-\sigma_{\rm t,\rm yield}^2}{2(\kappa_{\rm t}{+}\kappa_{_{\rm H}}{+}\kappa_{_{\rm H}}\tan^2\psiG)}\left(1{+}\tan^2\psiG\right)}
&\text{for }\  \arcsin\displaystyle{\frac{\sigma_{\rm t,\rm yield}}{\sqrt{2\kappa_{\rm t} a_{_{\rm I}}}}}\leq\psiG\leq \frac{\pi}{2}.
\end{cases}
\label{GcpsiG}
\end{align}
%
\REDD{By comparing the expression of $\alpha(\JMP{u})$ in \eqref{GcpsiG} with \eqref{alpha-split} and \eqref{a-mixity-dependence-}, and taking into account that $\alpha_{_{\rm ADHES}}=a_{_{\rm I}}$ is a constant, it is straightforward  to extract a somewhat cumbersome expression  of  $\alpha_{_{\rm PLAST}}(\JMP{u})$ for $\arcsin{\frac{\sigma_{\rm t,\rm yield}}{\sqrt{2\kappa_{\rm t} a_{_{\rm I}}}}}\leq\psiG$, whereas $\alpha_{_{\rm PLAST}}(\JMP{u})=0$ for $\psiG\leq\arcsin{\frac{\sigma_{\rm t,\rm yield}}{\sqrt{2\kappa_{t} a_{_{\rm I}}}}}$.}
{F}rom \eqref{GcpsiG} the  maximum value of $G_{\rm c}$  is given by, cf. \eqref{2-conditions},
\begin{equation}\label{a90}
a_{_{\rm II}}=G_{\rm c}(90^\circ)=a_{_{\rm I}}\left(1+\frac{\kappa_{\rm t}}{\kappa_{_{\rm H}}}\right) - \frac{\sigma^2_{\rm t,\rm yield}}{2\kappa_{_{\rm H}}}.
\end{equation}

In the LEBIM, we will test the functional dependence defined in \eqref{GcpsiG}, governed by
two parameters $\sigma_{\rm t,\rm{yield}}$ and $\kappa_{_{\rm H}}$, in addition to the parameters
$a_{_{\rm I}}$, $\kappa_{n}$ and $\kappa_{t}$ of a basic linear elastic-brittle model which is
(originally) insensitive to \REDD{fracture} mode mixity. Plots of the normalized fracture energy
$G_{\rm c}(\psiG)/a_{_{\rm I}}$ in Fig.~\ref{paramratio} qualitatively represent the behaviour
observed in experiments \cite{BanAsk00NFCI,HutSuo92MMCL,LieCha92ASIF,SwLiLo99ITAM,ERDC90FEBI}.

\begin{my-picture}{.9}{.41}{fig5:delam-mode-II}
\psfrag{MODE-I}{\footnotesize \frame{\,MODE I$^{^{^{}}}_{}$\,}}
\psfrag{MODE-II}{\footnotesize \frame{\,MODE II$^{^{^{}}}_{}$\,}}
\psfrag{[u]n}{\footnotesize $\JUMP{u}{\rm n}$}
\psfrag{[u]t}{\footnotesize $\JUMP{u}{\rm t}$}
\psfrag{stress}{\footnotesize $\sigma_{\rm n}$}
\psfrag{stress2}{\footnotesize $\sigma_{\rm t}$}
\psfrag{sd}{\footnotesize \begin{minipage}[t]{5em}
$\!\sqrt{2\kappa_{\rm n}^{}a_{_{\rm I}}}$
\\[-.2em]$=\hspace*{-.3em}{:}\,\sigma_{\mathrm{n,crit}}$\end{minipage}}
\psfrag{sd1}{\footnotesize \begin{minipage}[t]{5em}
$\!\!\sqrt{2\kappa_{\rm n}^{}a_{_{\rm I}}}$
\\[-.2em]$=\hspace*{-.3em}{:}\,\sigma_{\mathrm{n,crit}}$\end{minipage}}
\psfrag{sd2}{\footnotesize \begin{minipage}[t]{5em}
$\!\!\sqrt{2\kappa_{\rm t}^{}a_{_{\rm I}}}$
\\[-.2em]$=\hspace*{-.3em}{:}\,\sigma_{\mathrm{t,crit}}$\end{minipage}}
\psfrag{sd3}{\footnotesize $\hspace*{-1.7em}\sigma_{\mathrm{t,yield}}$}
\psfrag{s}{\footnotesize $\sqrt{2a_{_{\rm I}}/\kappa_{\rm n}^{}}$}
\psfrag{s2}{\small $\frac{\sigma_{\mathrm{t,yield}}}{\kappa_{\rm t}^{}}$}
\psfrag{s3}{\footnotesize $u_{_{\rm II}}$ from \eqref{ch5:delam-small-s-pi-II}}
\psfrag{s4}{\footnotesize  $\sqrt{2a_{_{\rm II}}/\kappa_{\rm t}^{}\!}$}
\psfrag{area}{\footnotesize area=$a_{_{\rm I}}$}
\psfrag{slope}{\footnotesize slope=$\kappa_{\rm n}^{}$}
\psfrag{slope3}{\footnotesize \hspace{1cm}slope=$\kappa_{\rm t}^{}$}
\psfrag{slope1}{\footnotesize slope=$2a_{_{\rm II}}/u_{_{\rm II}}^2$}
\psfrag{slope2}{\footnotesize{slope=}{\small$\frac{\kappa_{\rm t}\kappa_{_{\rm H}}}{\kappa_{\rm t}{+}\kappa_{_{\rm H}}}$}}
\psfrag{area2}{\footnotesize area=$a_{_{\rm II}}$}
\psfrag{pi0}{\footnotesize $\kappa_{_{\rm H}}\pi_{_{\rm II}}$}
\psfrag{evolution}{\footnotesize back-stress $\kappa_{_{\rm H}}\pi$}
\psfrag{of-pi}{\footnotesize of $\pi$}
\psfrag{LEBIM}{\footnotesize\bf LEBIM}
\psfrag{RIM}{\footnotesize\bf APRIM}
\hspace*{-.5em}\vspace*{-.1em}{\includegraphics[width=.89\textwidth]{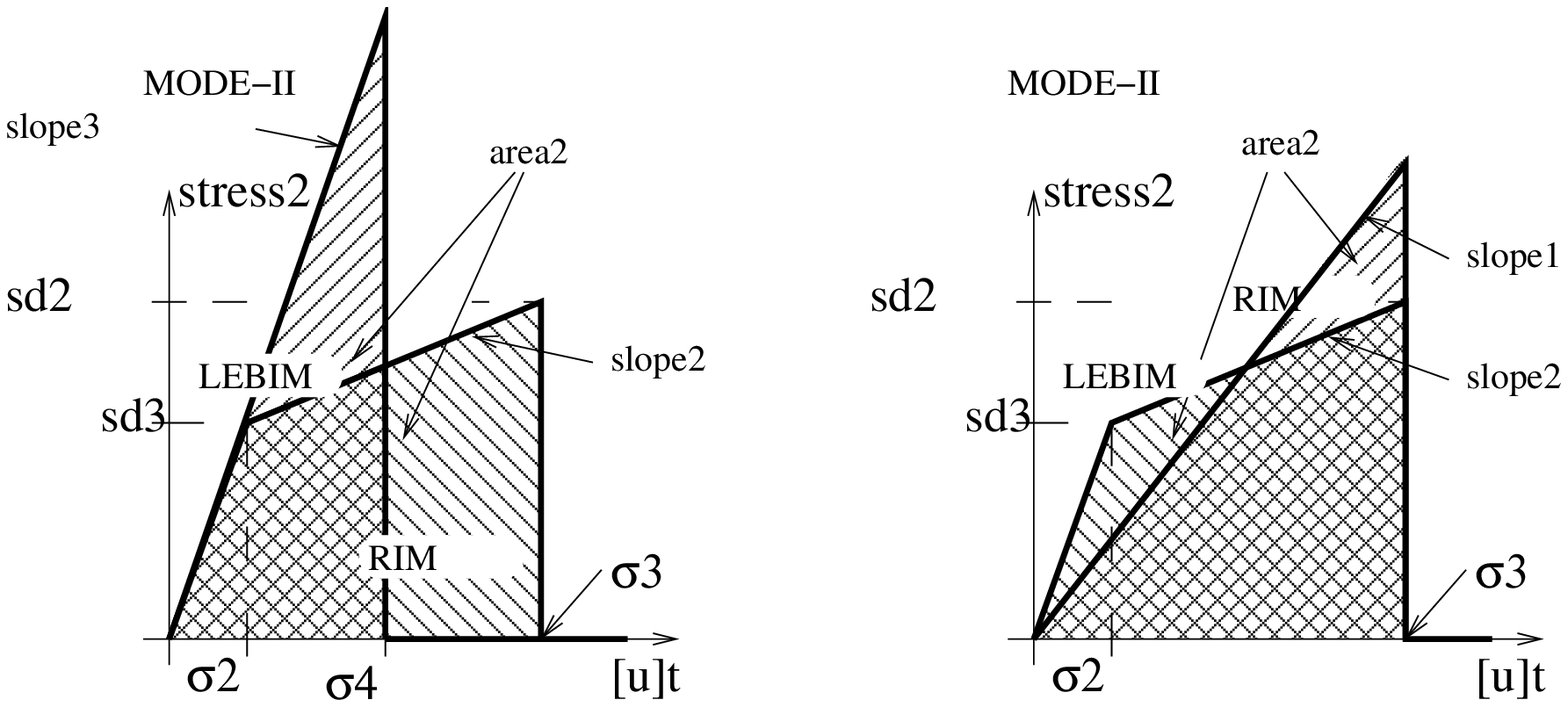}}
\end{my-picture}\nopagebreak
\begin{center}
{\small Fig.\,\ref{fig5:delam-mode-II}.\ }
\begin{minipage}[t]{.85\textwidth}\baselineskip=8pt
{\small
Two scenarios for LEBIM in Mode II after fitting  $\alpha$ by using
APRIM with \TTT a specific $\bbA$. \EEE
Left: $\bbA=
\TTT{\rm diag}(\kappa_{\rm n}^{},\kappa_{\rm t}^{})\EEE$.
Right: $\bbA={\rm diag}(\kappa_{\rm n}^{},2a_{_{\rm II}}/u_{_{\rm II}}^2)$
with $a_{_{\rm II}}$ and $u_{_{\rm II}}$ from \eqref{ch5:delam-small-a-s-II}.
}
\end{minipage}
\end{center}

Let us remind that, as mentioned in Sect.\,\ref{sec-model-A},
the above phenomenological relation \CHECK{\eqref{GcpsiG}}, motivated by an elasto-plastic behaviour of an adhesive layer, can   be rewritten in terms of $\psi_u$ or $\psi_{\sigma}$ instead of $\psiG$. For an actual interface of some stiffness $\kappa_{\rm t}^{}$, the dependence of the fracture energy on the  fracture-mode-mixity is controlled by two parameters: the yield stress $\sigma_\mathrm{t,yield}^{}$ and the hardening parameter $\kappa_{\scriptscriptstyle\textrm{H}}$. According to the plots presented in Fig.~\ref{paramratio}, the functional dependence of \CHECK{$\alpha(\psiG)(=G_c(\psiG))$} qualitatively represents the expected behaviour in view of the previous experimental results \cite{BanAsk00NFCI,HutSuo92MMCL,LieCha92ASIF,SwLiLo99ITAM,ERDC90FEBI}. It can easily be observed from these plots  that, as it was expected, $\sigma_\mathrm{t,yield}^{}$ has an influence on the threshold value of $\psiG$ for which $\alpha(\psiG)$ changes its behaviour from a constant function to an increasing function, while \REDD{both $\sigma_\mathrm{t,yield}^{}$ and} the hardening parameter $\kappa_{\scriptscriptstyle\textrm{H}}$  \REDD{control}  the  \CHECK{ value  of ratio $\alpha(\psiG)/a_{_{\rm I}}$,   in particular its maximum value according to \eqref{a90}}.

\begin{my-picture}{.99}{.37}{paramratio}
\hspace*{-1em}\vspace*{-.0em}{\includegraphics[width=.5\columnwidth]{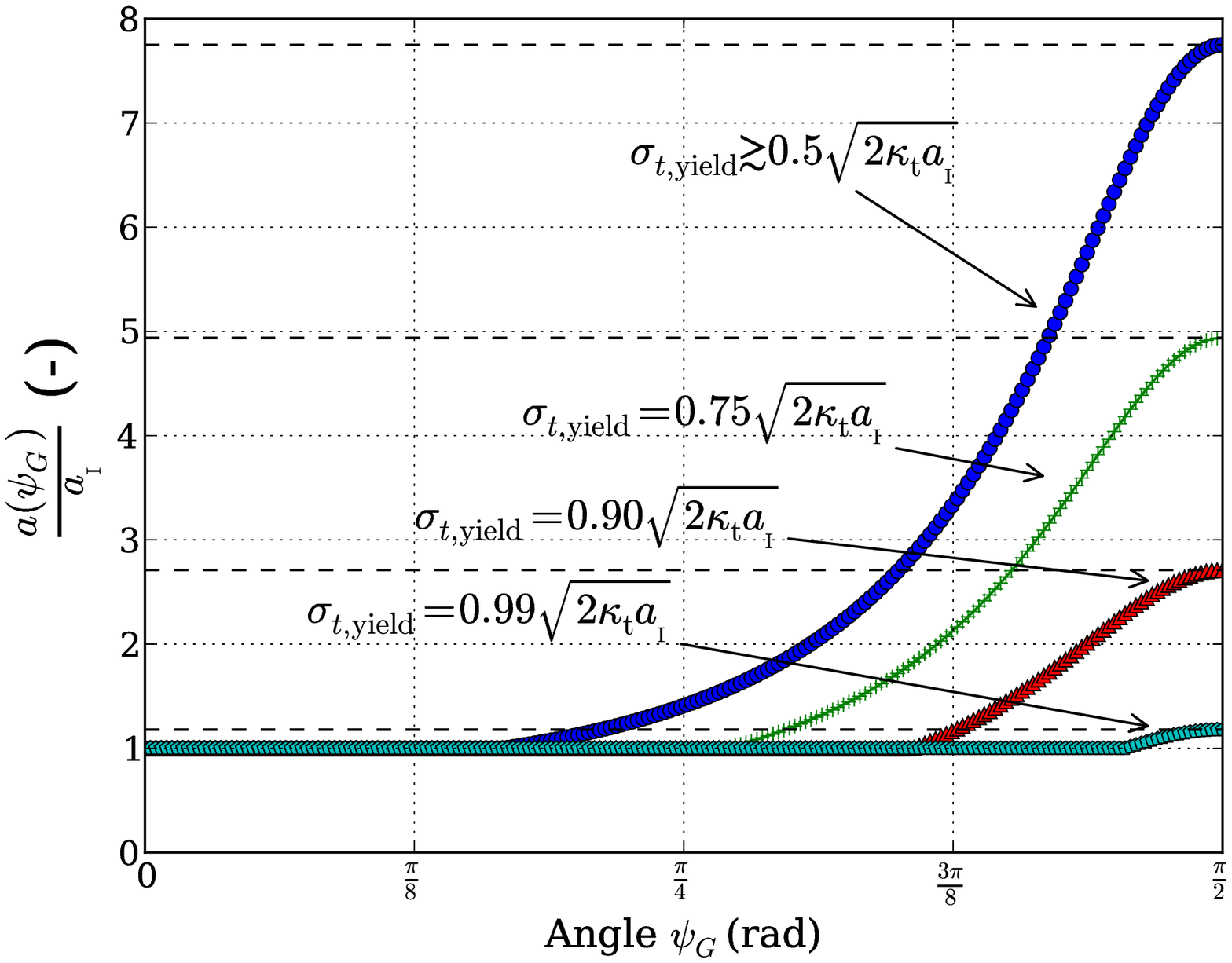}}{\includegraphics[width=.5\columnwidth]{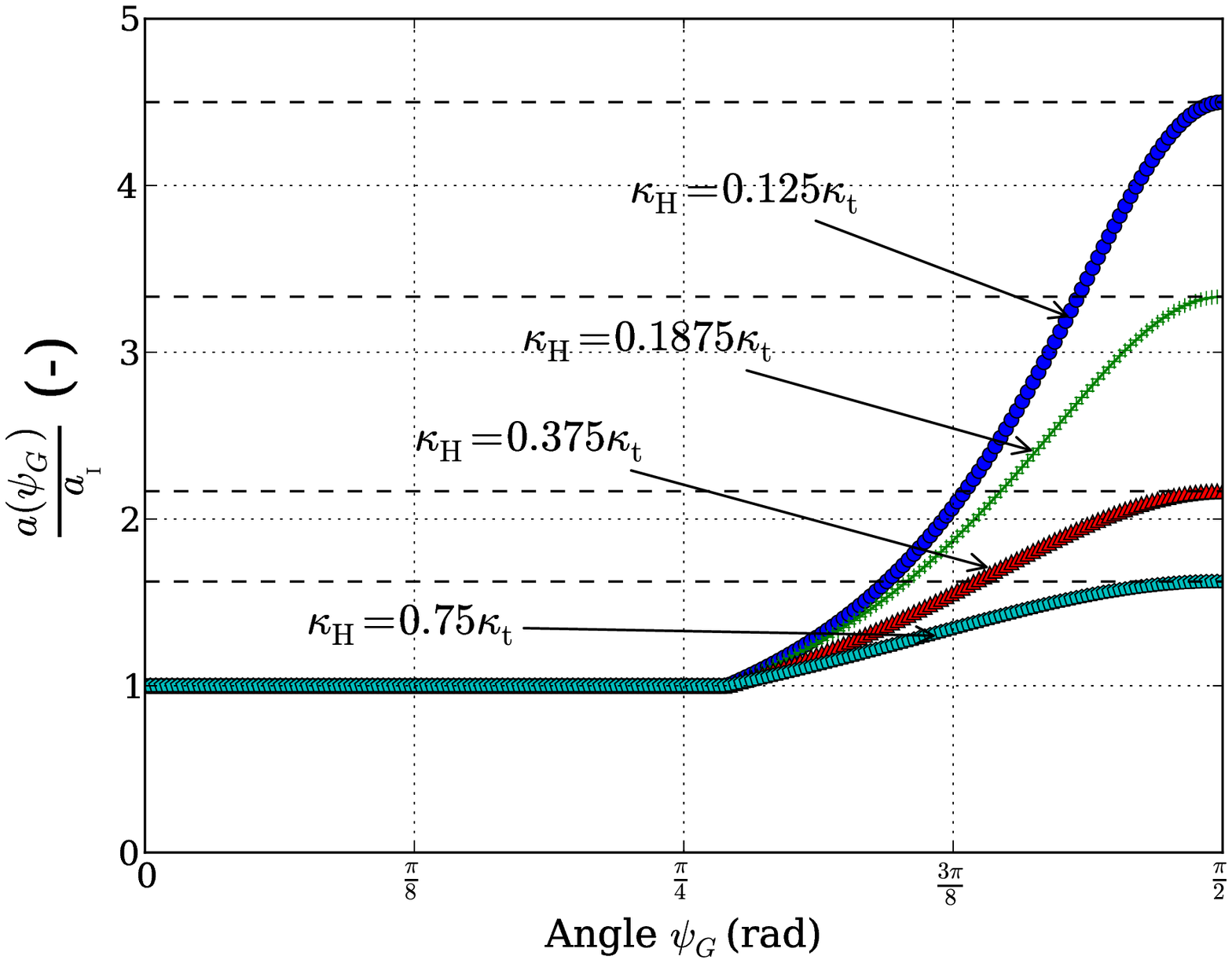}}
\end{my-picture}\begin{center}\hspace{-.3em}{\small Fig.\,\ref{paramratio}.\ }
\begin{minipage}[t]{
32em}\baselineskip=8pt
{\small
 Influence of $\sigma_\mathrm{t,yield}^{}$, considering
$\kappa_{\scriptscriptstyle\textrm{H}}{=}0.11\kappa_{\rm t}^{}$ (left), and of $\kappa_{\scriptscriptstyle\textrm{H}}$, considering $\sigma_\mathrm{t,yield}^{}{=}0.75\sqrt{2\kappa_{\rm t}^{}a_{_{\rm I}}}$ (right), on the
ratio $\alpha(\psi_G)/a_{_{\rm I}}$ given by \eqref{GcpsiG}.
In both plots, the respective $a_{_{\rm II}}/a_{_{\rm I}}$ value from \eqref{2-conditions} is plotted with the dashed line.}
\end{minipage}
\end{center}
\COMMENT{\LARGE I would recommend to change Left and Right to (a) and (b), which is much more efficient, clear and avoids confusion during typesetting in different journals.}

\section{Illustrative 2D simulations and comparison}\label{sec-simul}
\subsection*{Single\REDD{-}domain \REDD{example} with LEBIM and APRIM 
comparison}
\CHECK{A comparison of both models including varying fracture-mode-mixity of
delamination is shown} on a geometrically relatively simple but anyhow nontrivial
two-dimensional example
motivated by the pull-push shear experimental test used in engineering
practice \cite{CoCa11ES}. Intentionally, we use the same geometry,
shown in Fig.~\ref{fig_m1},
as in \cite{RoMaPa13QMMD,RoPaMa14LSAQ} in order to have a comparison of our
weak solution of the engineering non-associative visco-elastic model
with a maximally-dissipative
\COMMENT{\LARGE A NUMERICAL EVALUATION OF MAX-DISS DESIRED!!}
local solution
of the associative inviscid model presented
in \cite{RoPaMa14LSAQ}.
In contrast to the previous sections,
only one bulk domain is considered and $\GC$ is a part of its
boundary but this modification is straightforward;
alternatively, one may also think about
$\Omega_-$ as a completely rigid body in the previous setting.

\begin{my-picture}{.95}{.17}{fig_m1}
\psfrag{GN}{\footnotesize $\GN$}
\psfrag{GD}{\footnotesize $\GD$}
\psfrag{GC}{\footnotesize $\GC$}
\psfrag{elastic}{\footnotesize elastic body}
\psfrag{obstacle}{\footnotesize rigid obstacle}
\psfrag{adhesive}{\footnotesize adhesive}
\psfrag{LC}{$L_c$}
\psfrag{L}{\footnotesize $L=\ $250\,mm}
\psfrag{H}{\footnotesize $H=$}
\psfrag{12.5}{\scriptsize 12.5\,mm}
\psfrag{loading}{\footnotesize loading}
\hspace*{-.5em}\vspace*{-.1em}{\includegraphics[width=.95\textwidth]{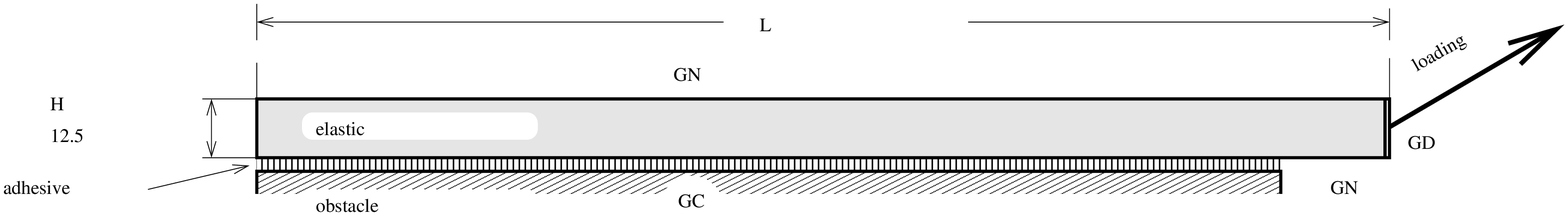}}
\end{my-picture}
\vspace*{-1.3em}
\begin{center}
{\small Fig.\,\ref{fig_m1}.\ }
\begin{minipage}[t]{.75\textwidth}\baselineskip=8pt
{\small
Geometry and boundary conditions of the problem considered.
The length of the initially glued part $\GC$ is $0.9L=225\,$mm,
the adhesive layer has zero thickness.}
\end{minipage}
\end{center}

Here $\Omega_+$ is a two-dimensional rectangular domain glued along the
\REPLACE{most of its bottom side $\GC$}{90\% of its bottom side $\GC$ to a rigid obstacle, } with the Dirichlet loading acting on the
right-hand side $\GD$ in the direction $(1,0.6)$, see  Fig.~\ref{fig_m1},
increasing linearly in time with velocity
$0.3\ $mm/s.

The bulk material is considered  linear, homogeneous, and isotropic
with \REDD{Young's} modulus $E=70$ GPa and Poisson's ratio $\nu=0.35$
(which corresponds to aluminum);
thus $\bbC_{ijkl}=\frac{\nu E}{(1{+}\nu)(1{-}2\nu)}\delta_{ij}\delta_{kl}
+\frac E{2(1{+}\nu)}(\delta_{ik}\delta_{jl}+\delta_{il}\delta_{jk})$
with $\delta_{ij}$ standing for the Kronecker symbol. \CHECK{In the LEBIM, the viscosity tensor
  $\bbD=\chi\bbC$ is considered} with a relaxation
time $\chi=0.001\,$s, which is very small in   relation to the
considered external-loading speed. Actually, we did not see any essential
difference for just merely elastic material with $\chi=0$, although the inviscid
variant of the LEBIM is not theoretically justified \REDD{even
as far as existence of solution concerns, and thus neither convergence
of any numerical scheme}.

For the adhesive, we took the normal stiffness $\kappa_{\rm n}^{}=$150 GPa/m, the
tangential stiffness with $\kappa_{\rm t}^{}=\kappa_{\rm n}^{}/2$, and the
Mode-I fracture toughness $a_{_{\rm I}}=187.5$ J/m$^2$. Furthermore,
the engineering LEBIM was fitted \DELETE{by \eqref{GCDEP2}} to  the APRIM using
$\kappa_{\scriptscriptstyle\textrm{H}}=\kappa_{\rm t}^{}/9$ and
$\sigma_\mathrm{t,yield}^{}\approx 0.79\sqrt{2\kappa_{\rm t}^{}a_{_{\rm I}}}$,
which satisfies inequalities \REDD{in} eq.~\eqref{2-sided-condition} and corresponds
to a rather moderate mode-sensitivity $\alpha(90^\circ)/a_{\rm I}\approx 4.36$.

\COMMENT{NOTHING SAID ABOUT BEM DISCRETIZATION AND TIME STEP ? Christos: not true see at the end of that section}


\vspace*{-14.5em}
\begin{my-picture}{.95}{.85}{fig_m3}
\begin{tabular}{ll}
\hspace*{2em}\LARGE $^{^{^{^{\mbox{\footnotesize $k$=50}}}}}$  & \hspace*{.5em}\vspace*{-.1em}$^{^{^{^{^{^{{\includegraphics[width=.6\textwidth]{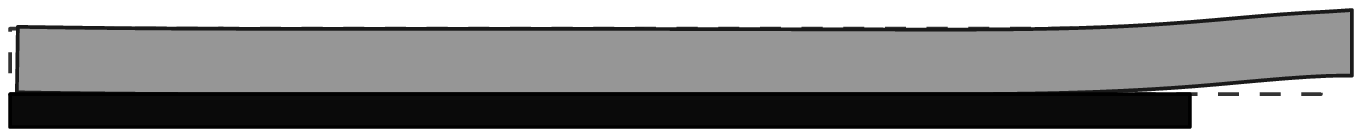}}}}}}}}$
\\[-1em]
\hspace*{2em}\LARGE $^{^{^{^{\mbox{\footnotesize $k$=87}}}}}$  & \hspace*{.5em}\vspace*{-.1em}$^{^{^{^{^{^{{\includegraphics[width=.6\textwidth]{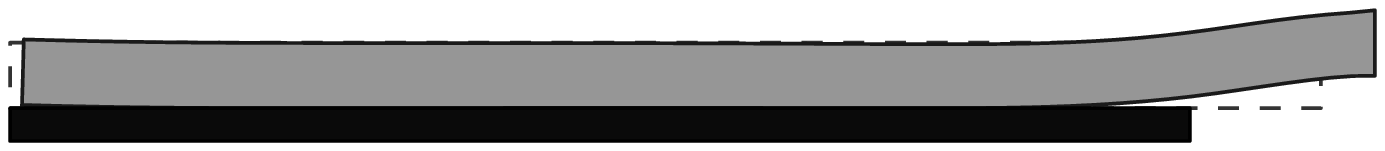}}}}}}}}$
\\[-1em]
\hspace*{2em}\LARGE $^{^{^{^{\mbox{\footnotesize $k$=107}}}}}$  & \hspace*{.5em}\vspace*{-.1em}$^{^{^{^{^{^{{\includegraphics[width=.6\textwidth]{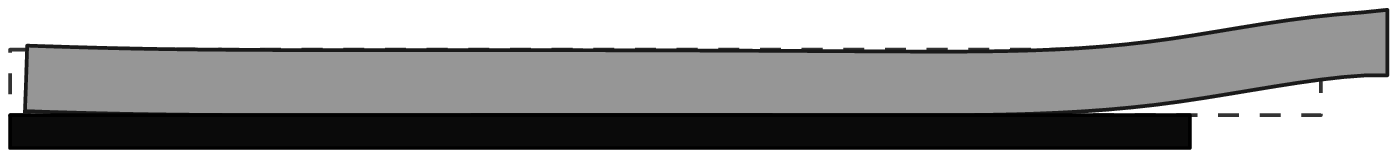}}}}}}}}$
\\[-1em]
\hspace*{2em}\LARGE $^{^{^{^{\mbox{\footnotesize $k$=163}}}}}$  & \hspace*{.5em}\vspace*{-.1em}$^{^{^{^{^{^{{\includegraphics[width=.6\textwidth]{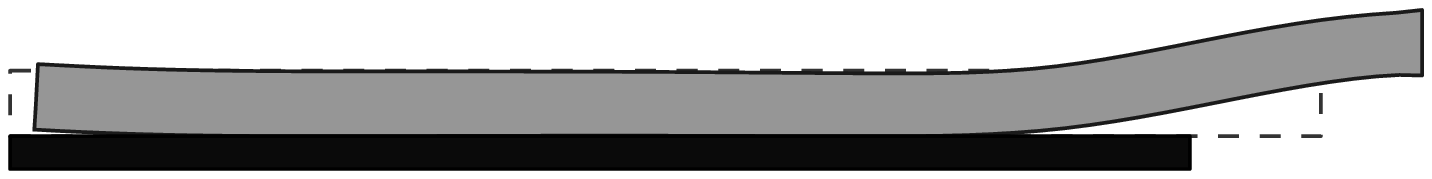}}}}}}}}$
\\[-1em]
\hspace*{2em}\LARGE $^{^{^{^{\mbox{\footnotesize $k$=198}}}}}$  & \hspace*{.5em}\vspace*{-.1em}$^{^{^{^{{\includegraphics[width=.6\textwidth]{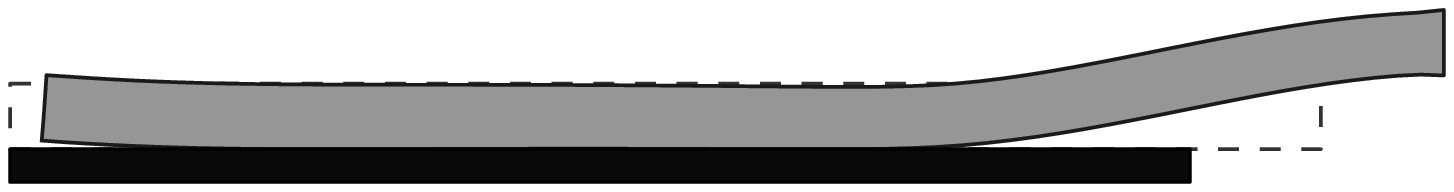}}}}}}$
\\[-.7em]
\hspace*{2em}\LARGE $^{^{^{^{\mbox{\footnotesize $k$=238}}}}}$  & \hspace*{.5em}\vspace*{-.1em}$^{^{^{^{{\includegraphics[width=.6\textwidth]{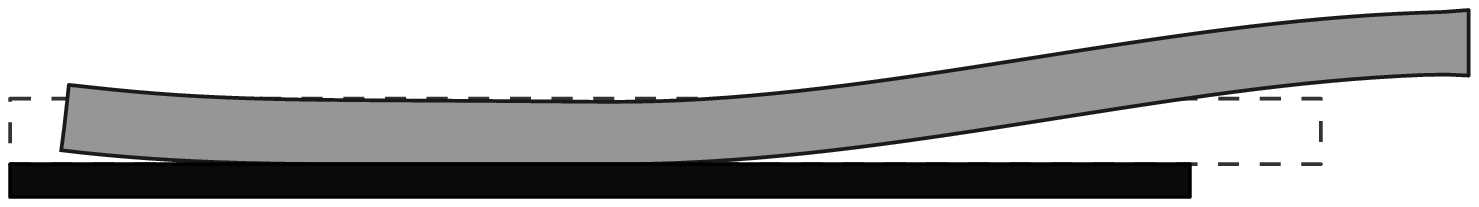}}}}}}$
\\[-.5em]
\hspace*{2em}\LARGE $^{^{^{^{\mbox{\footnotesize $k$=242}}}}}$  & \hspace*{.5em}\vspace*{-.1em}$^{^{^{^{{\includegraphics[width=.6\textwidth]{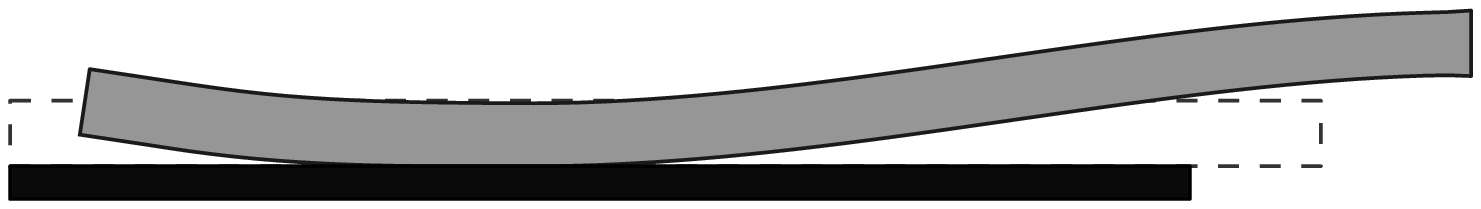}}}}}}$
\\[-.3em]
\hspace*{2em}\LARGE $^{^{^{^{\mbox{\footnotesize $k$=254}}}}}$  &
\hspace*{.5em}\vspace*{-.1em}$^{^{{\includegraphics[width=.6\textwidth]{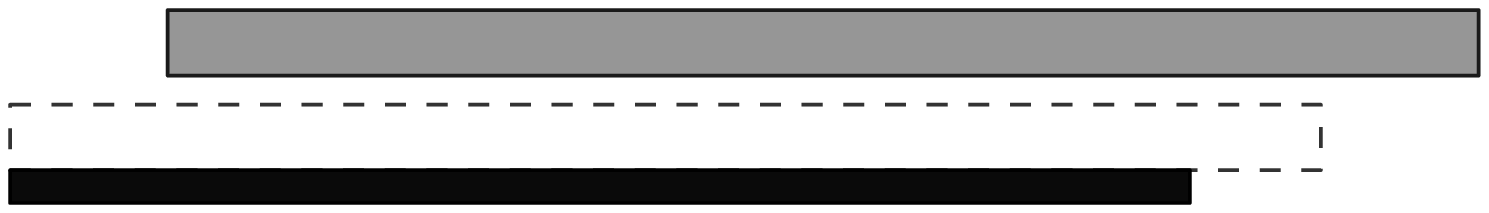}}}}$
\end{tabular}
\end{my-picture}
\vspace{17em}
\begin{center}
{\small Fig.\,\ref{fig_m3}.\,}
\begin{minipage}[t]{.89\textwidth}\baselineskip=8pt
{\small
Time evolution at eight snapshots of the deformed geometrical configuration
(solid lines) compared with the original configuration (dashed lines)
until complete delamination (displacement depicted magnified $100\,\times$).
Calculations performed by the algorithm LEBIM.
}
\end{minipage}
\end{center}

This example exhibits remarkably varying mode of delamination. At the
beginning, the delamination is performed by a mixed mode close to Mode I given
essentially by the direction of the Dirichlet loading, see Figure~\ref{fig_m1},
while later it turns rather to nearly pure Mode II. Yet, at the
very end of the process, due to elastic bending
the delamination starts performing also from the
left-hand side of the bar opposite to the loading side, and
thus again a mixed mode occurs. This relatively complicated mixed-mode
behaviour is depicted as a ``movie'' of 8 selected snapshots in
Figure~\ref{fig_m3} calculated by the LEBIM fitted to the
APRIM through \REPLACE{\eqref{GCDEP2}}{the model of Fig.~\ref{fig5:delam-mode-II}(Left)}, which would yield essentially similar picture.

The comparison of ``global'' quantities for the non-associative engineering
LEBIM with $\alpha(\cdot)$ from
 the models of Fig.~\ref{fig5:delam-mode-II}
with the associative (plasticity) APRIM is shown in Figures~\ref{fig_comp}
and \ref{fig_comp+}.
As it can be observed the two models, even if they are based on similar assumptions \CHECK{and their parameters are fitted to provide   similar responses}, they show quite different behaviours. The observed differences can be explained by the fact that in APRIM the stiffness of the  adhesive layer is progressively decreasing due to the presence of an interface plastic slip $\pi$, which may evolve at a portion of the interface  before the damage suddenly occurs. Nevertheless, in LEBIM there is no variation of the adhesive layer stiffnesses before the damage occurs. In particular,  this effect is clearly seen in Fig. \ref{fig_comp}, where the specimen stiffnesses begins to decrease with
respect to that in LEBIM once the plastic slip begins \TTT to \EEE evolve. It worths \TTT
mentioning \EEE
that the LEBIM mode-mixity distribution is in a good agreement with the expected  maximum  $\,a_{_{\rm II}}/a_{_{\rm I}}$ value.

Let us eventually briefly outline some implementation details \INSERT{for both models}.
We used the semi-implicit discretisation of the fractional-step type.
Using an equidistant partition of the time interval $[0,T]$ with a time step
$\tau>0$, for LEBIM, the discrete variant of \eqref{Biot-A} looks as
\begin{subequations}\label{Biot-A-disc}\begin{align}
&\partial_{\DT u}\mathcal{R}\Big(\frac{u_\tau^k{-}u_\tau^{k-1}}\tau\Big)
+\partial_u\mathcal{E}(k\tau,u_\tau^k,z_\tau^{k-1})\ni0
\ \ \text{ and }
\\&\partial_{\DT z}\mathcal{R}\Big(u_\tau^k;\frac{z_\tau^k{-}z_\tau^{k-1}}\tau\Big)
+\partial_z\mathcal{E}(u_\tau^k,z_\tau^k)\ni0,
\end{align}\end{subequations}
\COMMENT{\LARGE Is it ok to have $u_\tau^{k-1}$ in dissip. potential while $u_\tau^k$ in the stored energy?}
where $u_\tau^k$ denotes an approximation of $u(k\tau)$ and similarly for
$z_\tau^k$. This recursive scheme is to be solved for $k=1,2,...,T/\tau\in\N$.
Note that both inclusions in \eq{Biot-A-disc} are decoupled from
each other. For APRIM, the discrete variant of \eqref{Biot-B} looks as
\begin{subequations}\label{Biot-B-disc}\begin{align}\label{Biot-B-disc+}
&\partial_u\mathcal{E}(k\tau,u_\tau^k,z_\tau^{k-1},\pi_\tau^k)\!\ni\!0,
\ \ \text{ and }\ \
\partial_{\DT\pi}\mathcal{R}\Big(\frac{\pi_\tau^k{-}\pi_\tau^{k-1}}\tau\Big)
+\partial_\pi\mathcal{E}(u_\tau^k,z_\tau^{k-1},\pi_\tau^k)\!\ni\!0,
\ \text{ and }\ \ \
\\\label{Biot-B-disc++}
&\partial_{\DT z}\mathcal{R}\Big(\frac{z_\tau^k{-}z_\tau^{k-1}}\tau\Big)
+\partial_z\mathcal{E}(u_\tau^k,z_\tau^k,\pi_\tau^k)\!\ni\!0,
\end{align}\end{subequations}
Also this recursive scheme decouples \eqref{Biot-B-disc+} from
\eqref{Biot-B-disc++}.

\begin{my-picture}{.95}{.38}{fig_comp}
\hspace*{-1.em}\vspace*{-.1em}{\includegraphics[width=.5\textwidth]{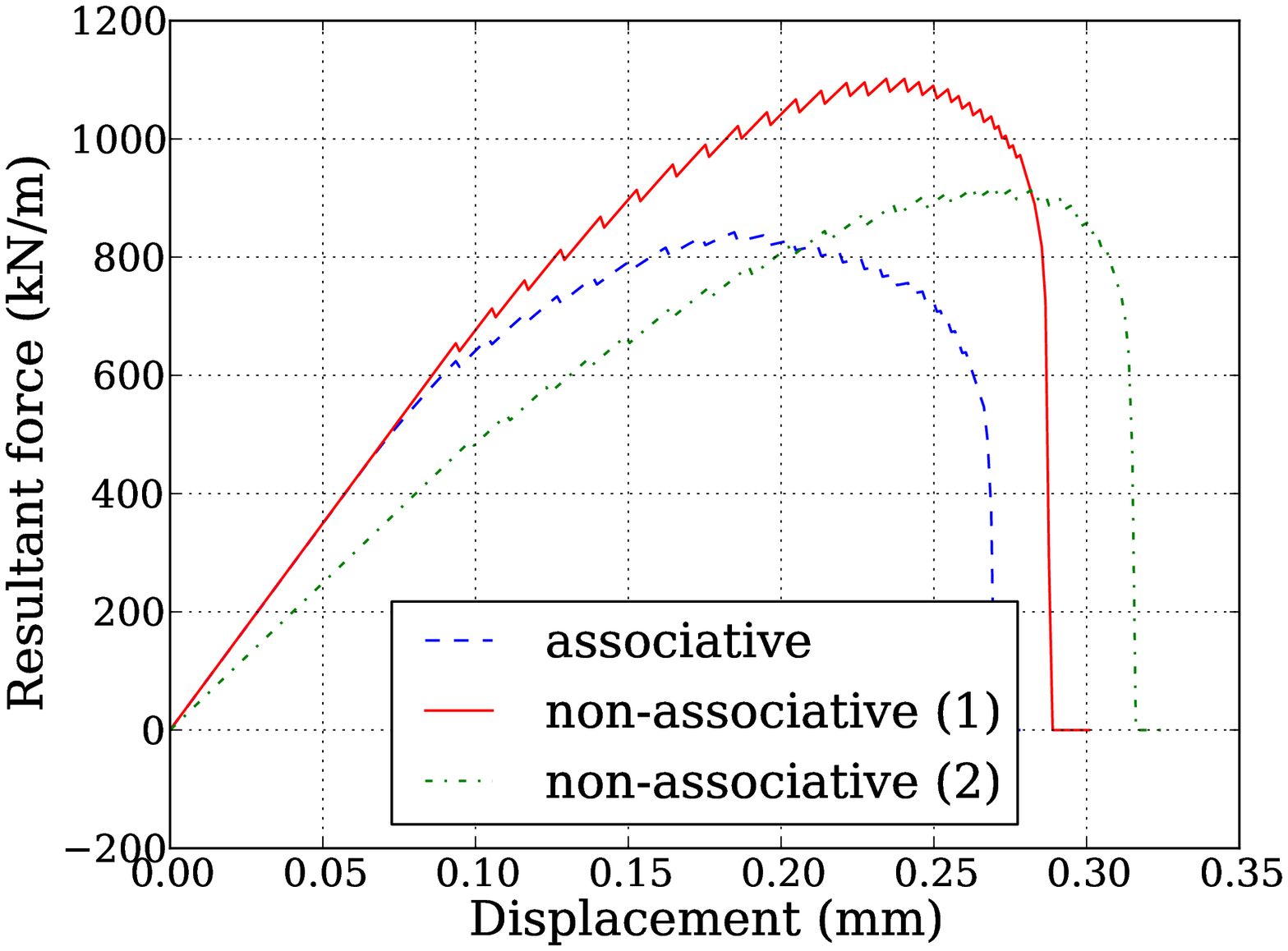}}
\hspace*{-.5em}\vspace*{-.1em}{\includegraphics[width=.5\textwidth]{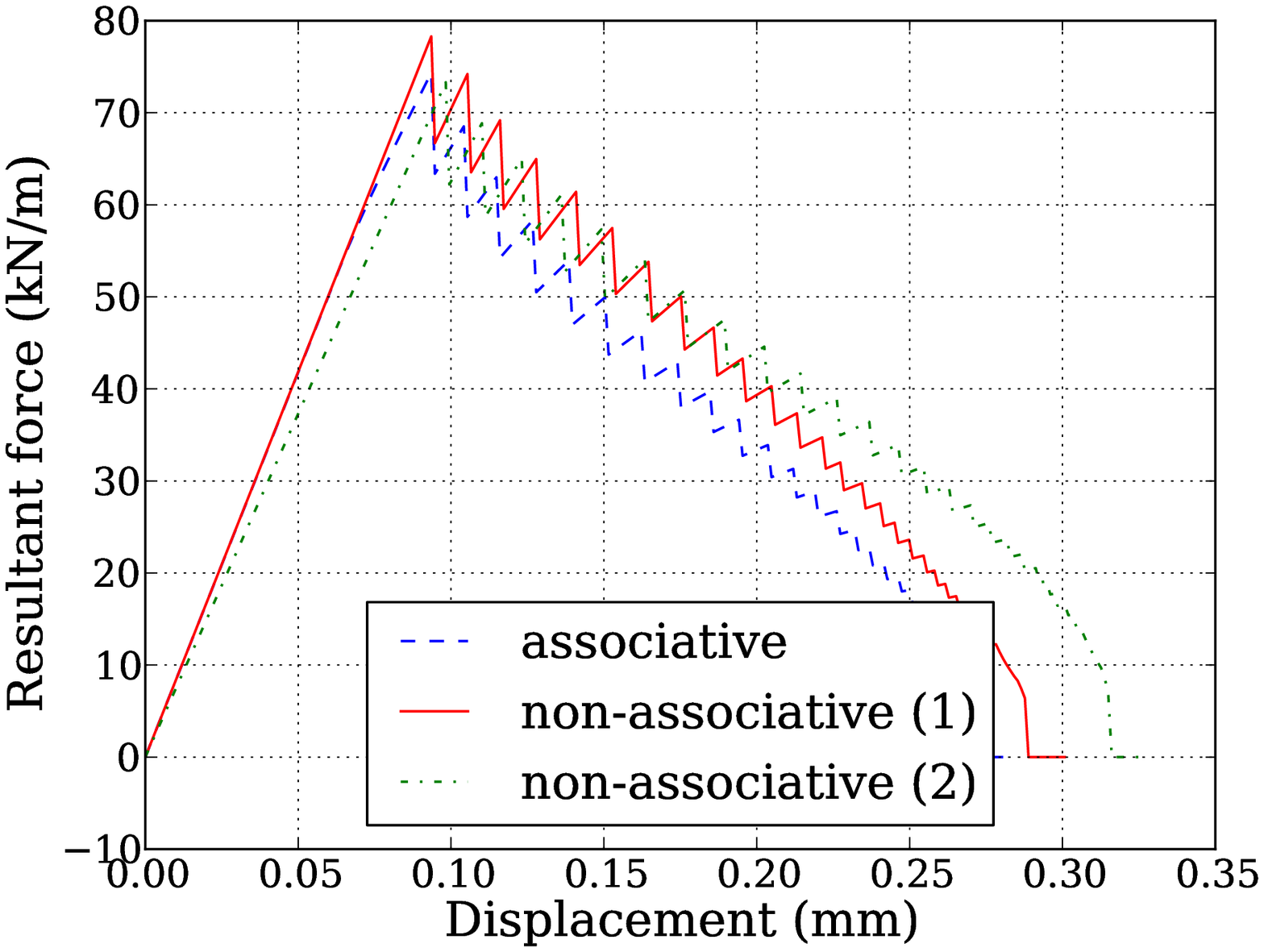}}
\end{my-picture}
\vspace{-1em}
\begin{center}
{\small Fig.\,\ref{fig_comp}.\ }
\begin{minipage}[t]{.8\textwidth}\baselineskip=8pt
{\small
The total force response evolving in time: the horizontal
(left) and the vertical (right) components.
}
\end{minipage}
\end{center}
\COMMENT{\LARGE Christos, I think, we have already discussed recently, with reference to differences in $a_{II}/a_{I}$ in Fig. 7(Right) because you evaluate in $a_{II}$ just dissipated energy not the hardening energy, which is the top triangle (clear grey color) in Fig. 1. Am I right? Could it be corrected? We possibly would obtain a better agreement in mode mixity.(Christos reply): Corrected.}
\COMMENT{\LARGE Another issue is the gap between the total and dissipated energy in APRIM in Fig. 7(left), at the end of their evolutions. Whereas, there is no such gap in the non-associative model. Why this difference? (Christos reply): To be honest I can not remember right now but I think it was something expected when we did the work for the Math. Models Methods Appl. Sci.\ldots}

\begin{my-picture}{.95}{.35}{fig_comp+}
\hspace*{-1.em}\vspace*{-.1em}{\includegraphics[width=.5\textwidth]{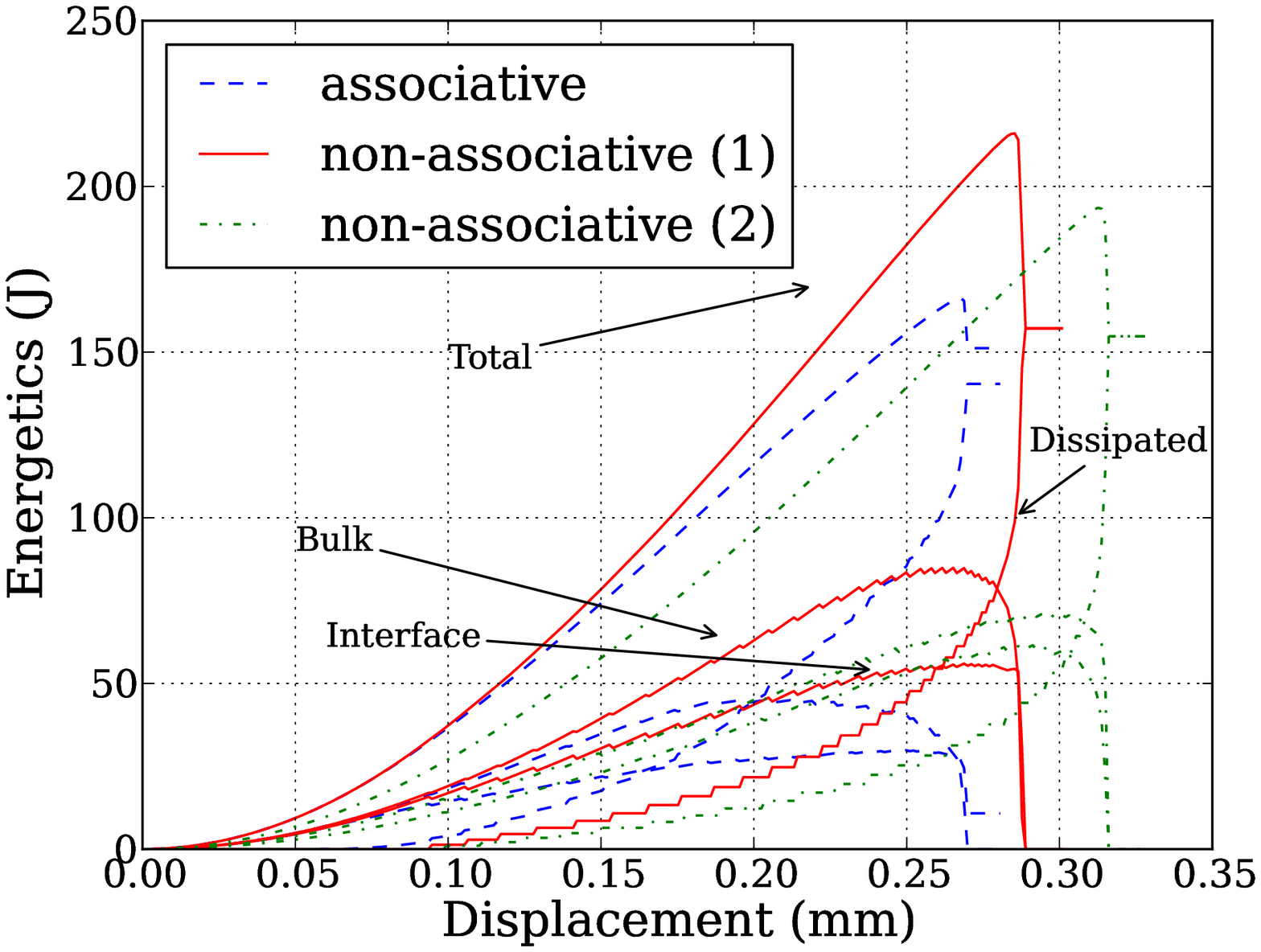}}
\hspace*{-.5em}\vspace*{-.1em}{\includegraphics[width=.5\textwidth]{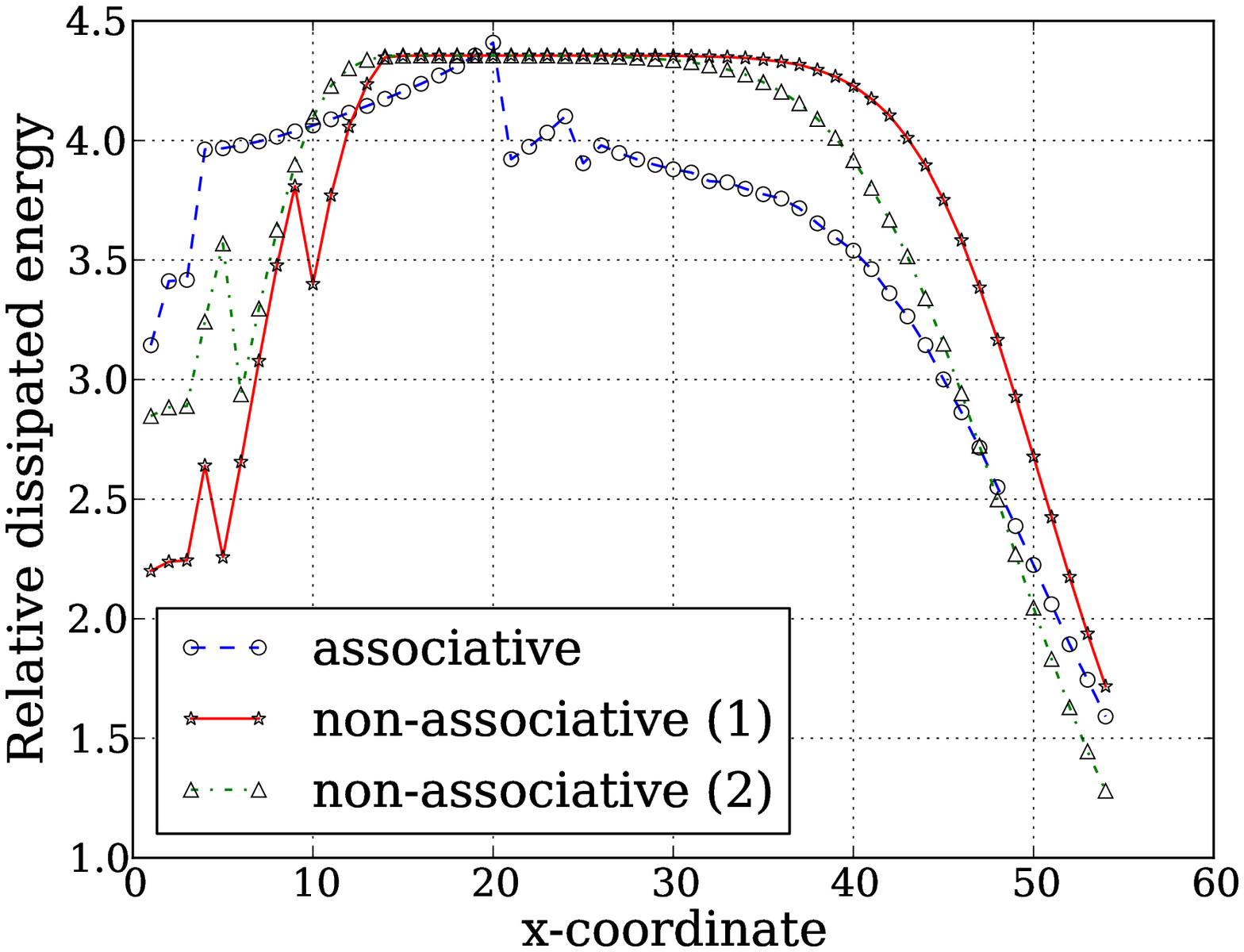}}
\end{my-picture}
\vspace{-1em}
\begin{center}
{\small Fig.\,\ref{fig_comp+}.\ }
\begin{minipage}[t]{
33em}\baselineskip=8pt
{\small
{\sf Left:} Time evolution of the energies:
the bulk and the interfacial parts of the stored energy
(together with the viscous contribution in the case of the LEBIM),
the interfacial dissipated energy due to delamination,
their sum\,=\,total energy.\COMMENT{OK??}
\\
{\sf Right:} Fracture-mode-mixity distribution along elements of $\GC$
evaluated \REDD{as the ratio of the overall dissipated energy to}
$a_{_{\rm I}}$ after the delamination has been completed:
value\,=\,1$\,\sim\,$Mode I, value\,=\,4.36\,=$\,a_{_{\rm II}}/a_{_{\rm I}}\!\sim\,$Mode II.

}
\end{minipage}
\end{center}

For details of spatial discretisation by P1/P0-boundary-element method
(P1 used for displacements, tractions and plastic slip $\pi$ discretizations,
whereas P0 used for the damage variable $z$  discretization) and the
outlined semi-implicit discretisation leading to a quadratic-programming
or even a linear-programming problems, as well as of computer implementation
in the case of both models LEBIM and APRIM, we refer to \cite{KruPaRo??MMS}
and \cite{RoMaPa13QMMD,RoPaMa14LSAQ}, respectively.
For BEM-implementation of the involved boundary-value problems
see also
\cite{TMGP11BEMA,RoPaMA13QACV,Mukh01BIEC,ParCan97BEMF}.

For the results presented in Figures~\ref{fig_m3}--\ref{fig_comp+},
we have used 54 elements on $\GC$, i.e.\ $h=4.1\bar{6}\,$mm
(=the size of a boundary segment in  the present uniform discretization), and
the time step $\tau = 5\,$ms for both LEBIM and APRIM simulations.

\begin{remark}[\em{Approximate maximum-dissipation principle}]\label{rem-AMDP}
\upshape
We already mentioned \REDD{at the end of Sect.\,\ref{sec-model-B}} that 
\eqref{Biot-B-disc} enjoys a guaranteed
stability and convergence towards weak (called also local) solutions
to the rate-independent problem \eqref{Biot-B}.
Yet, this class of
solutions may be very wide if the stored-energy functional
$\calE(t,\cdot,\cdot,\cdot)$ is nonconvex, which is our case here,
cf.\ \cite{Miel05ERIS} for a survey of variety of solution concepts,
some of them may exhibit nonphysically early ruptures due to
counting with a big energy cumulated in the elastic bulk or a nonphysical
tendency of sliding to less dissipative Mode I even in situations
where obviously Mode II should be preferred, cf.\ also
\cite{VoMaRo14EMDL}. It was discussed in
\cite{Roub??MDLS} that the sound Hill's maximum-dissipation principle
\cite{Hill48VPMP} can serve as an additional attribute which may
eliminate such unphysical weak solutions. In our case, it reads as
\begin{align}\label{e1:max-dis-principle}
\int_{\GC}\!\!\xi(t)\DT z(t)\,\d S
=\max_{\widetilde\xi\in K_z}\int_{\GC}\!\!\widetilde\xi\,\DT z(t)\,\d S
\ \ \ \text{ and }\ \ \ \int_{\GC}\!\!\zeta(t){\cdot}\DT\pi(t)\,\d S
=\max_{\widetilde\zeta\in K_\pi}\int_{\GC}\!\!\widetilde\zeta{\cdot}\DT\pi(t)\,\d S
\end{align}
for $t\in [0,T]$, where $\xi\in-\partial_z\calE(u(t),z(t),\pi(t))$ and
$\zeta\in-\partial_\pi\calE(u(t),z(t),\pi(t))$ are available driving forces
occurring already in \eqref{adhes-classic-APRIM} and $K_z:=\partial_{\dot{z}}\calR(0,0)$
and $K_\pi:=\partial_{\dot{\pi}}\calR(0,0)$ are the ``elasticity domains'' for
the internal variable $z$ and $\pi$, i.e.\ here
\COMMENT{\LARGE Vlado: I do not understand how $K_{\pi}$ is fixed (see eq. (5.4)), as in kinematic hardening the ``elastic region'' is changing with plasticity evolution; it does not change its size but translates. Surely I am missing something}
\begin{align}
K_z=\big\{\xi\!\in\!L^{1}(\GC);\ \xi\ge -a_1 \ \text{ on }\GC\big\}\ \text{ and }\
K_\pi=\big\{\zeta\!\in\!L^{\infty}(\GC;\R^{d-1});\ |\zeta|\le\sigma_{\mathrm{t,yield}^{}}
\ \text{ on }\GC\big\}.
\end{align}
This principle holds rather automatically for a.a.\ \TTT (=\,almost all) time
instants \EEE $t$ for any weak solution
to \eqref{adhes-classic-APRIM} which has $\DT z$ and $\DT\pi$ absolutely
continuous and then the integrals in \eqref{e1:max-dis-principle}
have the standard Lebesgue sense. On the other hand, during sudden ruptures
where $\DT z$ or $\DT\pi$ may concentrate in space and time to be only
a general measure, the validity of \eqref{e1:max-dis-principle}
is not automatic and even analytically the integrals in
\eqref{e1:max-dis-principle} looses a sense even as dualities. It has
been proposed in \cite{Roub??MDLS} to consider \eqref{e1:max-dis-principle}
rather integrated over the time interval $[0,T]$, which then reads,
when also re-writing ``max'' over the elasticity domains in
\eqref{e1:max-dis-principle} equivalently as total variations, as
\begin{align}\label{e1:max-dis-principle+}
\int_{\overlineSC}\!\!\xi(t)\,\d z(t)\d S=\int_{\GC}\!\!\REDD{a_1}\big(z(T){-}z_0\big)\,\d S
\ \ \ \text{ and }\ \ \ \int_{\overlineSC}\!\!\zeta(t){\cdot}\,\d\pi(t)\,\d S
=\sigma_\mathrm{t,yield}^{}\!\int_{\overlineSC}\!\!\big|\DT\pi(t)\big|\,\d S\d t
\end{align}
where the notations $\int\cdot\,\d z(t)$ \REDD{and 
$\int\cdot\,\d\pi(t)$} stand for \TTT (lower)
Pollard-Moore-Stieltjes integrals as used in \cite{MieRou15RIST},
\REDD{which are}
a Pollard-Moore modification of \EEE the so-called lower
Riemann-Stieltjes integral, cf.\ e.g.\ \cite{Rudi76PMA}, while
the last integral \REDD{is am integral of} 
a total variation of the measure $\DT\pi$.
Here, as suggested in \cite{Roub??MDLS}, we apply \eqref{e1:max-dis-principle+}
to the left-continuous piecewise-constant approximate solutions
obtained by \REPLACE{\eqref{Biot-A-disc}}{\eqref{Biot-B-disc}},
let us denote them by
$\bar u_\tau$, $\bar z_\tau$, and  $\bar\pi_\tau$. Such a maximum dissipation
principle for the approximate solution is expected to hold only approximately,
if at all. The lower \TTT Pollard-Moore\EEE-Stieltjes integral can then be explicitly
evaluated as well as the right-hand sides of \eqref{e1:max-dis-principle+},
which yields
\begin{subequations}\label{e1:AMDP}\begin{align}
&\int_{\SC}\!\!\bar\xi_\tau(t)\,\d\bar z_\tau(t)\d S
=\sum_{k=1}^{T/\tau}\int_{\GC}\!\!\xi_\tau^{k-1}(z_\tau^k-z_\tau^{k-1})\,\d S
\ \ {\stackrel{\mbox{\large\bf ?}}{\sim}}\ \
\int_{\GC}\!\!\CHECK{a_1}\big(z_\tau^{T/\tau}{-}z_0\big)\,\d S
\ \ \ \text{ and}
\\&\nonumber\int_{\SC}\!\!\bar\zeta_\tau(t)\,\d\bar\pi_\tau(t)\d S
=\sum_{k=1}^{T/\tau}\int_{\GC}\!\!\zeta_\tau^{k-1}(\pi_\tau^k-\pi_\tau^{k-1})\,\d S
\ \ {\stackrel{\mbox{\large\bf ?}}{\sim}}\ \
\sigma_\mathrm{t,yield}^{}\!\sum_{k=1}^{T/\tau}\int_{\GC}\!\!
\big|\pi_\tau^k-\pi_\tau^{k-1}\big|\,\d S
\\&\hspace{20.7em}=
\sigma_\mathrm{t,yield}^{} \!\int_{\SC}\!\!\big|\DT{\bar\pi}_\tau(t)\big|\,\d S\d t
\end{align}\end{subequations}
where $\xi_\tau^k\in-\partial_z\calE(u_\tau^k,z_\tau^{k-1},\pi_\tau^k)$ and
$\zeta_\tau^k\in-\partial_\pi\calE(u_\tau^k,z_\tau^k,\pi_\tau^k)$,
and where ''${\stackrel{\mbox{\bf ?}}{\sim}}$'' means that the
equality holds at most only asymptotically for $\tau\to \INSERT{0}$ but
even this is rather only desirable and not always valid.
See \cite{RoPaMa14LSAQ,VoMaRo14EMDL}, for such sort of
a-posteriori justification of the approximate local solution obtained
by fractional-step discretisation as a stress-driven solution.
\COMMENT{HERE SOME FIGURE WITH NUMERICAL RESULTS evaluating the
residua (possibly spatially or time distributed in
\eqref{e1:AMDP} WOULD BE NICE}
\COMMENT{Christos: we may use the plots from M3AS article with appropriate
references or elsewhere only references}
\end{remark}



\subsection*{Multi\REDD{-}domain example with LEBIM model \COMMENT{TO BE COMPLETED \ldots (?)}}\MARGINOTE{Here I give a rough mesh alternative while I will replace it with a finer one the following weeks. In the worst scenario I will replace it after the reviewing process! In any case I think it would be nice to include such an ``advanced'' example, in order to show efficiency of LEBIM.}
The influence of the mixed-mode behavior, as well as, the efficiency of the LEBIM   approximation, will be illustrated on an
example of the Mixed-Mode Flexure (MMF) test \cite{VC06EFM,tr+mk+jz,TMSGP13CM} shown in Fig.~\ref{mmfig}.
The specimen of the MMF test consists of two
aluminium arms, bonded with a layer of resin adhesive.
Detailed specimen dimensions and boundary conditions realizing the loading
are shown in Figure~\ref{mmfig}; \CHECK{it should be noted the point support
is on the bordeline of compatibility (and, in fact, not compatible) with
traces of $H^1$-functions in dimension 2 but anyhow the numerical simulations
did not see this slight discrepancy}. The material properties of the bulk and interface
are the same as in the previous example, \REDD{however considering $\sigma_\mathrm{t,yield}^{}\approx 0.56\sqrt{2\kappa_{\rm t}^{}a_{_{\rm I}}}$ and inviscid conditions}.
Plane strain conditions \REDD{are} considered for the problem, while a uniform mesh of  188 linear boundary elements is used to model the geometry. \COMMENT{SOMETHING SHOULD BE SAID ABOUT THE TIME STEP AND LOADING VELOCITY, AS WE INCLUDE NUMBERS OF SNAPHSOTS. THERE IS A PORBLEM WITH FIGURE NUMBERING HERE, CAUSED THAT PREVIOUSLY WE HAVE NOT MOVED THE LATEX-COUNTER OF FIGURES.(Christos answer): At time step $k=500$ imposed disp 2.5 mm. Visocity very small as in the first example, that is relaxation time 0.001 s.}
\TTT The evolution of the delamination of this double-stripe specimen in shown in Figure~\ref{fig_def2}. \EEE

\begin{my-picture}{.95}{.2}{mmfig}
\psfrag{GN}{\footnotesize $\GN$}
\psfrag{GD}{\footnotesize $\GD$}
\psfrag{GC}{\footnotesize $\GC$}
\psfrag{elastic+}{\scriptsize elastic body $\Omega_+$}
\psfrag{elastic-}{\scriptsize elastic body $\Omega_-$}
\psfrag{obstacle}{\footnotesize rigid obstacle}
\psfrag{adhesive}{\footnotesize adhesive}
\psfrag{L}{\footnotesize $120\,$mm}
\psfrag{3}{\scriptsize 3\,mm}
\psfrag{2}{\scriptsize 2\,mm}
\psfrag{60}{\scriptsize 60\,mm}
\psfrag{37}{\scriptsize 37\,mm}
\psfrag{10}{\scriptsize 10\,mm}
\psfrag{loading}{\footnotesize loading}
\vspace*{-3.3em}\includegraphics[width=.88\textwidth]{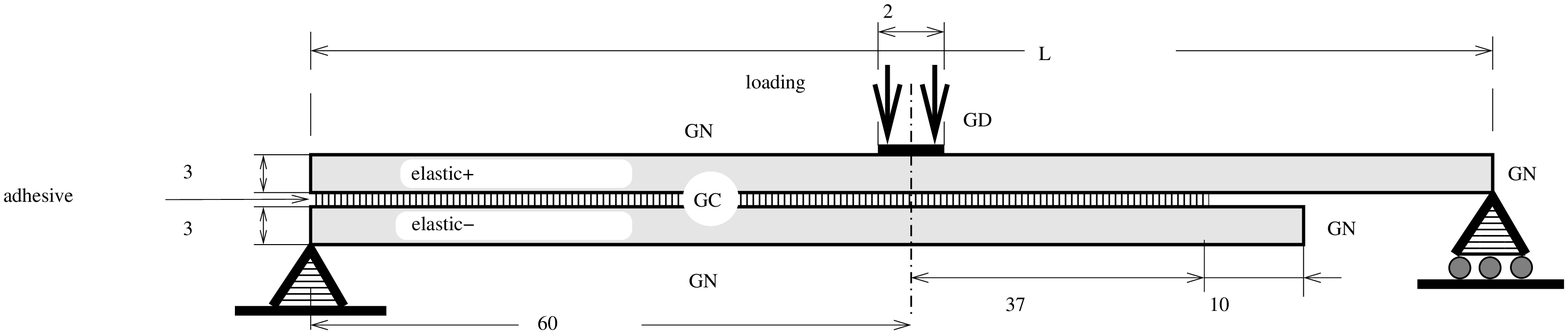}
\end{my-picture}
\vspace{-1em}
\begin{center}
{\small Fig.\,\ref{mmfig}.\,}
\begin{minipage}[t]{25em}\baselineskip=8pt
{\small
Specimen configuration in the Mixed-Mode Flexure test.
}
\end{minipage}
\end{center}

\bigskip\bigskip


\vspace*{-19em}
\begin{my-picture}{.95}{.95}{fig_def2}
\begin{tabular}{ll}
\hspace*{2em}\LARGE $^{^{^{^{\mbox{\footnotesize $k$=215}}}}}$  & \hspace*{.5em}\vspace*{-.1em}$^{^{^{^{^{^{{\includegraphics[width=.66\textwidth]{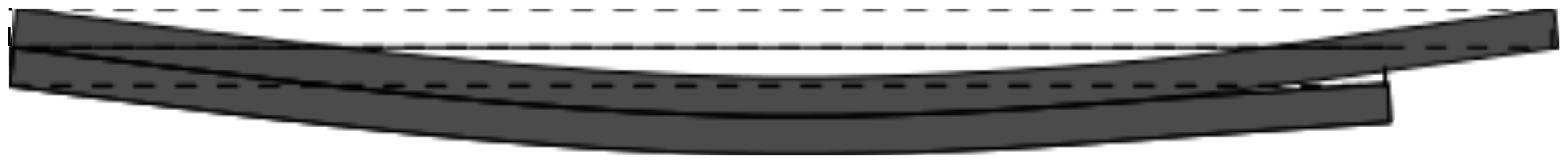}}}}}}}}$
\\[-1em]
\hspace*{2em}\LARGE $^{^{^{^{\mbox{\footnotesize $k$=217}}}}}$  & \hspace*{.5em}\vspace*{-.1em}$^{^{^{^{^{^{{\includegraphics[width=.66\textwidth]{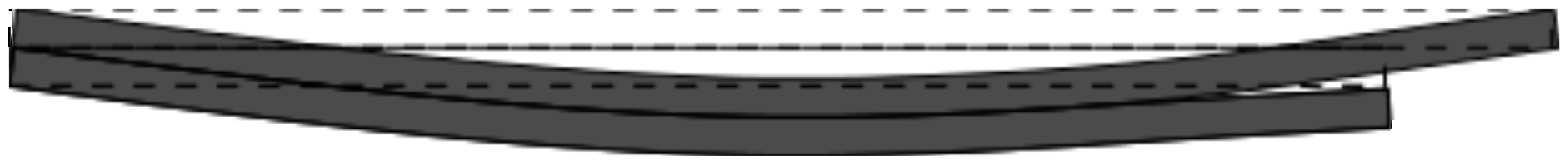}}}}}}}}$
\\[-1em]
\hspace*{2em}\LARGE $^{^{^{^{\mbox{\footnotesize $k$=219}}}}}$  & \hspace*{.5em}\vspace*{-.1em}$^{^{^{^{^{^{{\includegraphics[width=.66\textwidth]{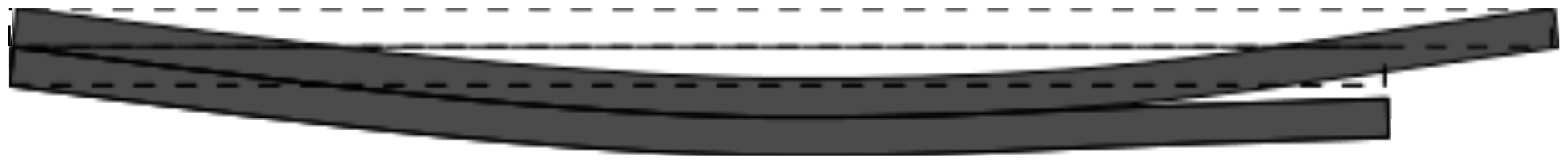}}}}}}}}$
\\[-1em]
\hspace*{2em}\LARGE $^{^{^{^{\mbox{\footnotesize $k$=221}}}}}$  & \hspace*{.5em}\vspace*{-.1em}$^{^{^{^{^{^{{\includegraphics[width=.66\textwidth]{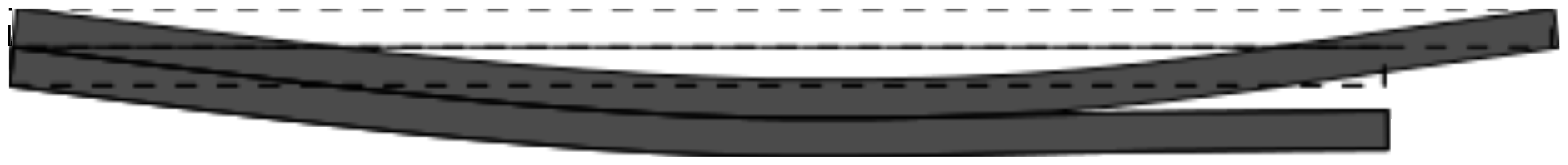}}}}}}}}$
\\[-1em]
\hspace*{2em}\LARGE $^{^{^{^{\mbox{\footnotesize $k$=223}}}}}$  & \hspace*{.5em}\vspace*{-.1em}$^{^{^{^{{\includegraphics[width=.66\textwidth]{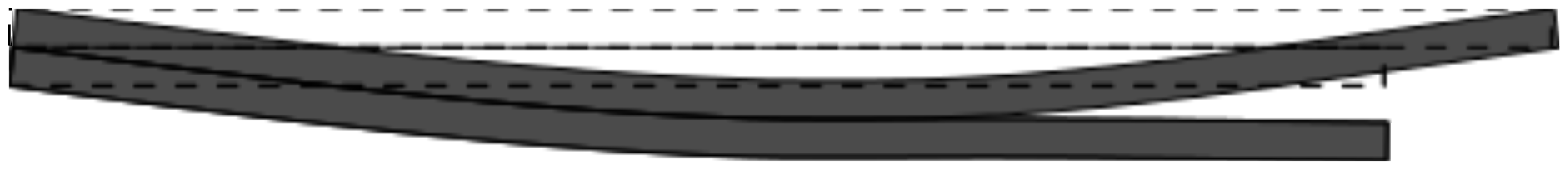}}}}}}$
\\[-.7em]
\hspace*{2em}\LARGE $^{^{^{^{\mbox{\footnotesize $k$=225}}}}}$  & \hspace*{.5em}\vspace*{-.1em}$^{^{^{^{{\includegraphics[width=.66\textwidth]{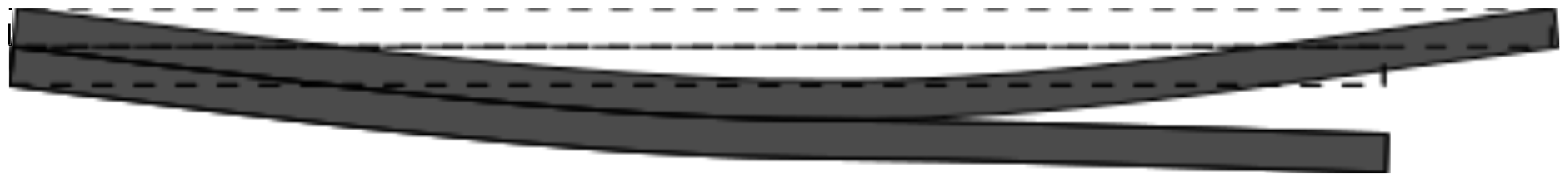}}}}}}$
\\[-.5em]
\hspace*{2em}\LARGE $^{^{^{^{\mbox{\footnotesize $k$=227}}}}}$  & \hspace*{.5em}\vspace*{-.1em}$^{^{^{^{{\includegraphics[width=.66\textwidth]{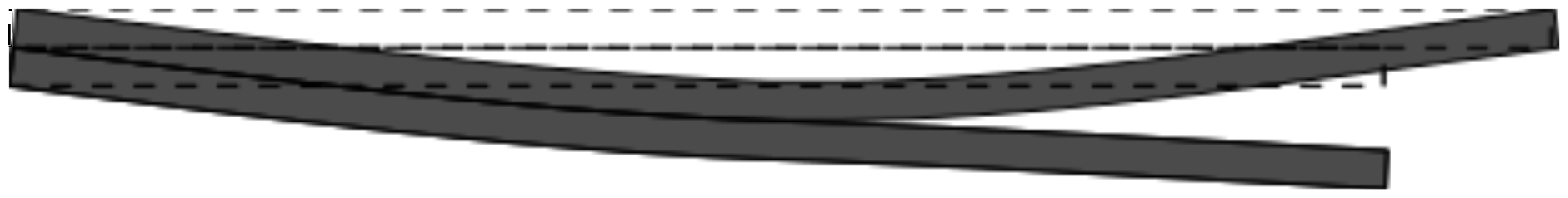}}}}}}$
\\[-.3em]
\hspace*{2em}\LARGE $^{^{^{^{\mbox{\footnotesize $k$=229}}}}}$  &
\hspace*{.5em}\vspace*{-.1em}$^{^{{\includegraphics[width=.66\textwidth]{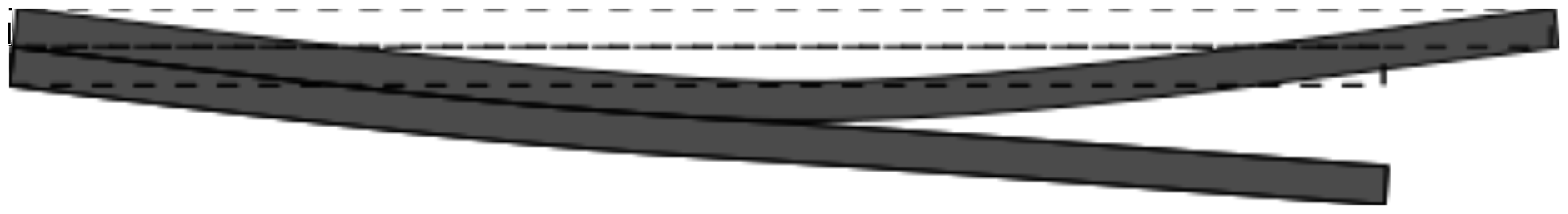}}}}$
\\[-.3em]
\hspace*{2em}\LARGE $^{^{^{^{\mbox{\footnotesize $k$=231}}}}}$  &
\hspace*{.5em}\vspace*{-.1em}$^{^{{\includegraphics[width=.66\textwidth]{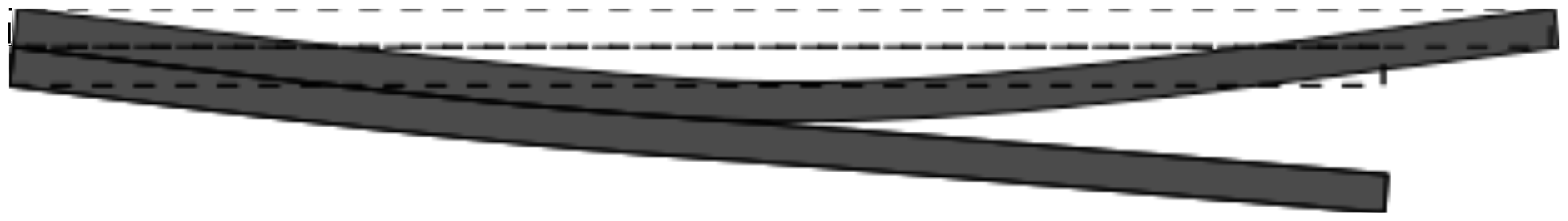}}}}$
\end{tabular}
\end{my-picture}
\vspace{17em}
\begin{center}
{\small Fig.\,\ref{fig_def2}.\,}
\begin{minipage}[t]{.89\textwidth}\baselineskip=8pt
{\small
Time evolution at ten snapshots of the deformed geometrical configuration
(solid lines) compared with the original configuration (dashed lines, displacement
depicted magnified $5\,\times$). Calculations performed by the algorithm of LEBIM.
}
\end{minipage}
\end{center}

Finally, both scenarios for the LEBIM model shown in Fig.~\ref{fig5:delam-mode-II}
are compared. The behaviour of their solutions  is shown in Fig.~\ref{fig_ex2_load},
where the energetics of the delamination evolution are summarized (left) together
with the load-deflection curves (right). \TTT Comparing this behaviour with the former
\REDD{test} from Fig.~\ref{fig_m3}, whose results are depicted \REDD{in}
Figs.~\ref{fig_comp} and \ref{fig_comp+}, we can see that  the scenario (1) \REDD{yields again} a bit earlier rupture than the scenario (2), \REDD{as could be expected from the stress-relative displacements laws of these scenarios shown in Fig.~\ref{fig5:delam-mode-II}. The length of the fractured part of the adhesive layer is the same in both scenarios, and the total dissipated energy is only slightly larger in the scenario (1) due to a slightly larger fracture mode mixity at the beginning of the crack propagation.} 
\EEE

\begin{my-picture}{.95}{.38}{fig_ex2_load}
\hspace*{-1.em}\vspace*{-.1em}{\includegraphics[width=.5\textwidth]{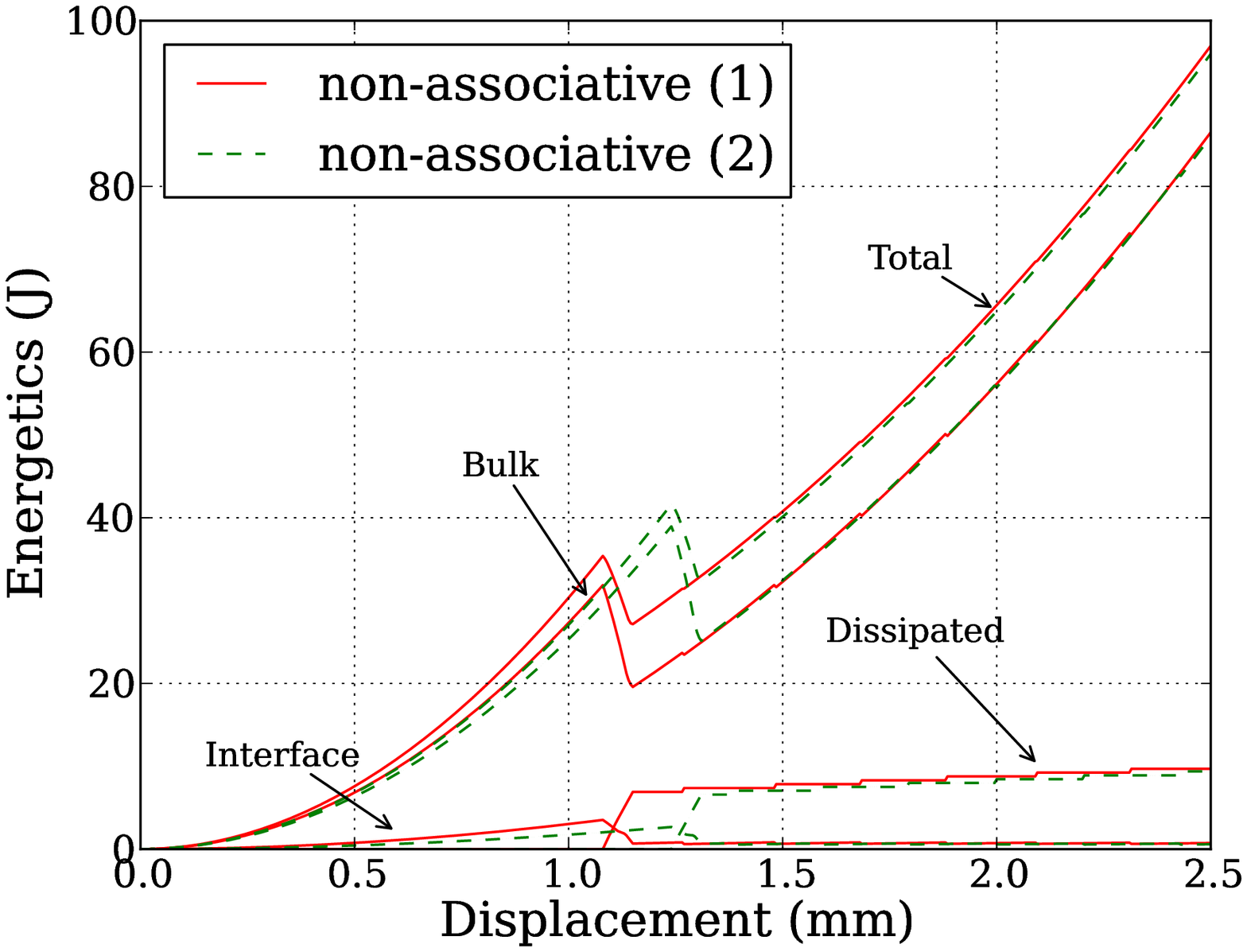}}
\hspace*{-.5em}\vspace*{-.1em}{\includegraphics[width=.5\textwidth]{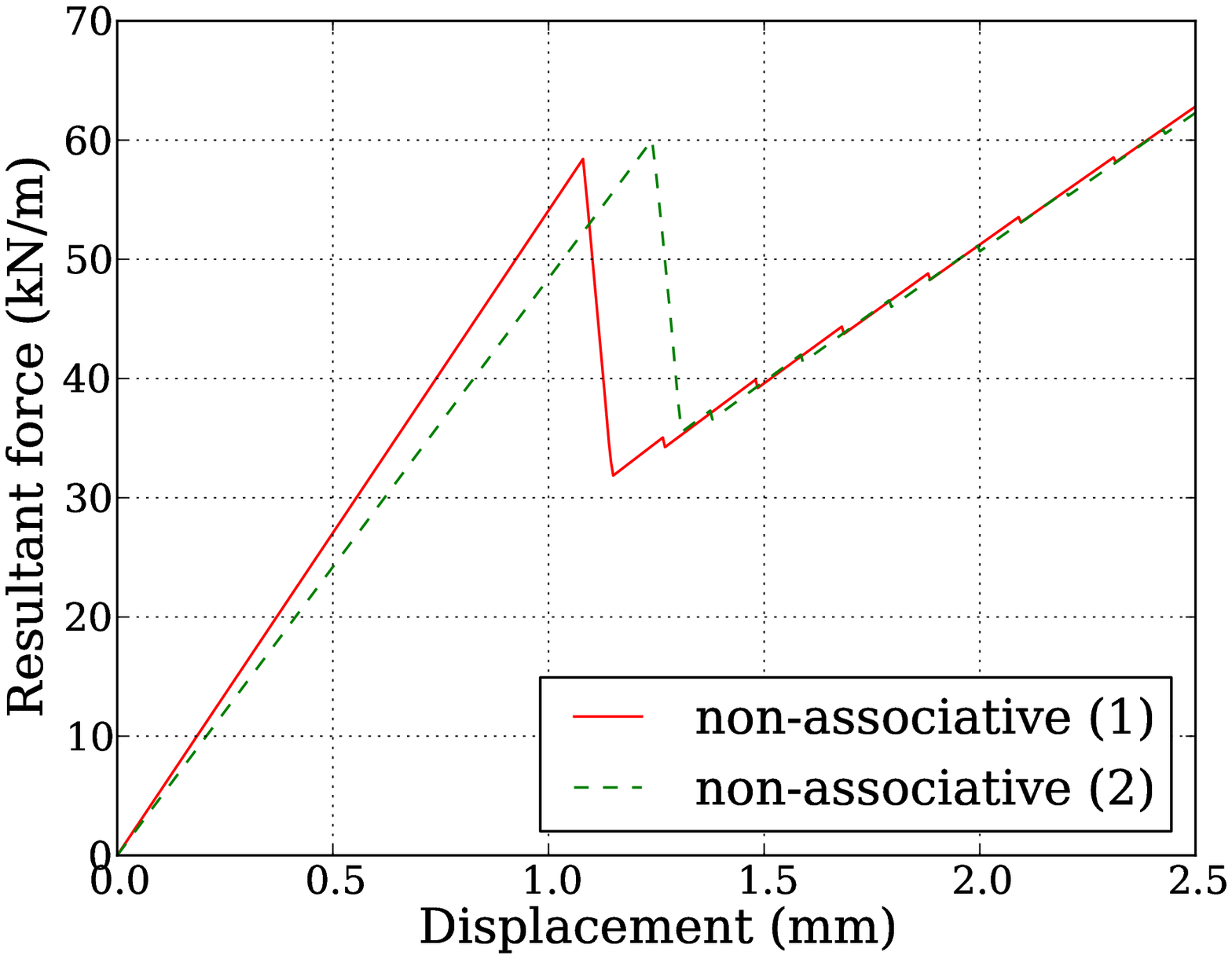}}
\end{my-picture}
\vspace{-1em}
\begin{center}
{\small Fig.\,\ref{fig_ex2_load}.\ }
\begin{minipage}[t]{.8\textwidth}\baselineskip=8pt
{\small
Comparison of the evolution of energies (left) and  load-deflection curves (right) for
both scenarios of Fig.~\ref{fig5:delam-mode-II}.
}
\end{minipage}
\end{center}

\subsection*{Acknowledgments}
This research has been supported by the Spanish Ministry of Economy and
Competitiveness and European Regional Development Fund (Project MAT2012-37387), the Junta de Andaluc\'{\i}a and
European Social Fund (Project of Excellence P08-TEP-04051),
and from the Czech Science Foundation through the grants 201/10/0357,
201/12/0671,
and 13-18652S,
together with the institutional support RVO: 61388998 (\v CR).

\bigskip\bigskip

\bibliographystyle{elsarticle-num}
\baselineskip=13pt
\noindent\textbf{\large References}
\bibliography{theBib}

\end{document}